\documentclass[11pt,a4paper]{article}
\usepackage[english]{babel}
\usepackage{amsmath,amsfonts,amssymb}
\usepackage{epsfig}
\usepackage{theorem}
\usepackage{subfigure}
\usepackage{pst-plot}
\usepackage{multirow}
\usepackage{titlesec}
\usepackage{indentfirst}
\usepackage{pst-tree}
\usepackage{lscape}

\theoremstyle{break}

\begin{document}

\begin{center}{\textbf{\Large A comparison of formulations and solution
methods for the Minimum-Envy Location Problem. Additional results}}
\end{center}

\long\def\symbolfootnote[#1]#2{\begingroup%
\def\thefootnote{\fnsymbol{footnote}}\footnote[#1]{#2}\endgroup}

\begin{center}
Inmaculada Espejo$^1$, Alfredo Mar\'{\i}n$^{2}$, Justo Puerto$^3$,
Antonio Rodr\'{\i}guez-Ch\'{\i}a$^1$
\end{center}

\begin{center} \small
$^1$ Departamento de Estad\'{\i}stica e Investigaci\'on Operativa, Universidad de C\'adiz, Spain,\\
$^2$ Departamento de Estad\'{\i}stica e Investigaci\'on Operativa, Universidad de Murcia, Spain,\\
$^3$ Departamento de Estad\'{\i}stica e Investigaci\'on Operativa, Universidad de Sevilla, Spain,\\
\end{center}

\symbolfootnote[0]{\textit{E-mails:} inmaculada.espejo@uca.es, amarin@um.es,
puerto@us.es, antonio.rodriguezchia@uca.es}

\section{Formulations}
\label{sec:form}

\begin{eqnarray}
\nonumber\label{F1} \hbox{(F1)} & \min & \sum_{i=1}^{M-1} \sum_{j=i+1}^M e_{ij} \\
\nonumber  & \hbox{s.t.}& e_{ij}\ge z_i-z_j, \quad \forall i=1,\ldots,M-1,j=i+1,\ldots,M \\
\nonumber  & & e_{ij}\ge z_j-z_i, \quad \forall i=1,\ldots,M-1,j=i+1,\ldots,M \\
\nonumber  & & z_i + \sum_{\ell :\ O_{i\ell }\le k-1}(k-O_{i\ell })y_{\ell } \ge k, \quad \forall i\in A,k=1,\ldots ,M-p+1 \\
\nonumber  & & z_i + (M-p+1-O_{ik})y_k \le M-p+1, \quad \forall i,k\in A:\ O_{ik}\le M-p \\
\nonumber  & & \sum_{j=1}^M y_j = p \\
\nonumber  & & y_j \in \{0,1\}, \quad \forall j\in A.
\end{eqnarray}

\begin{eqnarray}
\nonumber\label{F2} \hbox{(F2)} & \min & \sum_{i=1}^M (2i-M-1) \sum_{k=2}^{M-p+1} x_{ik} \\
\nonumber  & \hbox{s.t.} & z_{ik} \ge z_{i,k+1}, \quad \forall i\in A,k=2,\ldots,M-p \\
\nonumber  & & z_{iO_{ij}}-z_{i,O_{ij}+1} \le y_j, \quad \forall i,j\in A:\ 2\le O_{ij} \le M-p \\
\nonumber  & & 1-z_{i2} \le y_j, \quad \forall i,j\in A:\ O_{ij}=1 \\
\nonumber  & & z_{i,M-p+1} \le y_j, \quad \forall i,j\in A:\ O_{ij}=M-p+1 \\
\nonumber  & & z_{i,O_{ij}+1} + y_j \le 1, \quad \forall i,j\in A:\ O_{ij} \le M-p \\
\nonumber  & & \sum_{i=1}^M x_{ik} = \sum_{i=1}^M z_{ik}, \quad \forall k=2,\ldots,M-p+1 \\
\nonumber  & & x_{ik} \ge x_{i-1,k}, \quad \forall i=2,\ldots,M,k=2,\ldots,M-p+1 \\
\nonumber  & & \sum_{j=1}^M y_j=p \\
\nonumber  & & y_j \in \{0,1\}, \quad \forall j\in A \\
\nonumber  & & x_{ik}\in \{0,1\}, \quad \forall i\in A,\ k=2,\ldots ,M-p+1.
\end{eqnarray}

\begin{eqnarray}
\nonumber\label{F3} \hbox{(F3)} & \min &  \sum_{i=1}^M \sum_{q=1}^{M-1}2d_{iq}+\sum_{q=1}^{M-1} 2 q t_q - (M-1) \sum_{i=1}^M z_i \\
\nonumber  & \hbox{s.t.} &z_i + \sum_{\ell:\ O_{i\ell}\le k-1}(k-O_{i\ell})y_\ell \ge k, \quad \forall i\in A,k=1,\ldots ,M-p+1 \\
\nonumber  & & z_i + (M-p+1-O_{ik}) y_k \le M-p+1, \quad \forall i,k\in A:\ O_{ik}\le M-p+1 \\
\nonumber  & & \sum_{j=1}^M y_j = p \\
\nonumber  & & d_{iq} \ge -t_q + z_i, \quad \forall q=1,\ldots,M-1,\ i\in A \\
\nonumber  & & d_{iq} \ge 0 \quad \forall q=1,\ldots ,M-1,i\in A \\
\nonumber  & & y_j \in \{0,1\} ,\quad \forall j\in A.
\end{eqnarray}

\begin{eqnarray}
\nonumber\label{F4} \hbox{(F4)} & \min & \sum_{i=1}^M
\sum_{q=1}^{M-1} 2 d_{iq} + \sum_{q=1}^{M-1} 2 q t_q - (M-1) \sum_{i=1}^M \sum_{k=2}^{M-p+1} z_{ik} + M(1-M) \\
\nonumber  & \hbox{s.t.} & z_{ik} \ge z_{i,k+1}, \quad \forall i\in A, k=2,\ldots,M-p \\
\nonumber  & & z_{iO_{ij}}-z_{i,O_{ij}+1} \le y_j, \quad \forall i,j\in A:\ 2\le O_{ij} \le M-p \\
\nonumber  & & 1-z_{i2} \le y_j, \quad \forall i,j\in A:\ O_{ij}=1 \\
\nonumber  & & z_{i,M-p+1} \le y_j, \quad \forall i,j\in A:\ O_{ij}=M-p+1 \\
\nonumber  & & \sum_{j=1}^M y_j = p \\
\nonumber  & & z_{i,O_{ij}+1} + y_j \le 1, \quad \forall i,j\in A:\ O_{ij} \le M-p \\
\nonumber  & & d_{iq} \ge -t_q + \sum_{k=2}^{M-p+1} z_{ik} + 1, \quad \forall q=1,\ldots,M-1, \;i\in A \\
\nonumber  & & d_{iq} \ge 0, \quad \forall q=1,\ldots,M-1, \;i\in A \\
\nonumber  & & y_j \in \{0,1\}, \quad \forall j\in A.
\end{eqnarray}

\begin{eqnarray}
\nonumber\label{F5.1} \hbox{(F5.1)} & \min & \sum_{i=1}^M (2i-M-1) x_i \\
\nonumber & \hbox{s.t.} & z_i + \sum_{\ell:\ O_{i\ell}\le k-1}(k-O_{i\ell})y_\ell \ge k, \quad \forall i\in A,k=1,\ldots ,M-p+1 \\
\nonumber & & z_i + (M-p+1-O_{ik})y_k \le M-p+1, \quad \forall i,k\in A:\ O_{ik}\le M-p+1 \\
\nonumber & & \sum_{j=1}^M y_j=p \\
\nonumber & & x_i \le x_{i+1}, \quad \forall i=1,\ldots,M-1\\
\nonumber & & \sum_{i=1}^M x_i=\sum_{i=1}^M z_i \\
\nonumber & & \sum_{i=k}^M x_i \ge \sum_{i \in S} z_i, \quad \forall k=2,\ldots,M,\ \forall S\subset A:\ |S|=M-k+1\\
\nonumber & & y_j \in \{0,1\} , \quad \forall j\in A.
\end{eqnarray}

\begin{eqnarray}
\nonumber\label{F5.2}  \hbox{(F5.2)} & \min & \sum_{i=2}^M 2\theta_i - (M-1)\sum_{i=1}^M z_i \\
\nonumber  & \hbox{s.t.} & z_i + \sum_{\ell:\ O_{i\ell} \le k-1}(k-O_{i\ell}) y_\ell \ge k,
\quad \forall i\in A,k=1,\ldots, M-p+1 \\
\nonumber  & & z_i+ (M-p+1-O_{ik})y_k \le M-p+1, \quad \forall i,k\in A:\ O_{ik}\le M-p+1 \\
\nonumber  & & \sum_{j=1}^M y_j=p \\
\nonumber  & & \theta_k \ge \sum_{i\in S} z_i, \quad \forall k=2,\ldots ,M,\ \forall S\subset A:\  |S|=M-k+1\\
\nonumber  & & y_j \in \{0,1\} ,\quad \forall j \in A.
\end{eqnarray}

\section{Computational Results}
\label{sec:cresult}

The previous formulations have been compared by means of a
computational study (see Tables \ref{freeSS}, \ref{nofreeSS} and
\ref{random}).

\begin{table}
$$
\hspace{-2cm}
\begin{array}{|c|crrrr|rrrrr|rrrrr|}
\hline \multicolumn{1}{|c}{} & \multicolumn{5}{|c|}{M=20} &
\multicolumn{5}{|c|}{M=30} &
\multicolumn{5}{|c|}{M=40} \\
\hline
  & p & \multicolumn{1}{c}{LP} & \multicolumn{1}{c}{\bar{t}} & \multicolumn{1}{c}{\sigma_t} &
\multicolumn{1}{c}{n} & \multicolumn{1}{|c}{p} & \multicolumn{1}{c}{LP} & \multicolumn{1}{c}{\bar{t}} & \multicolumn{1}{c}{\sigma_t} & \multicolumn{1}{c}{n}  & \multicolumn{1}{|c}{p} & \multicolumn{1}{c}{LP} & \multicolumn{1}{c}{\bar{t}} & \multicolumn{1}{c}{\sigma_t} & \multicolumn{1}{c|}{n} \\
\hline
\hbox{(F1)}   &2 &0.0 &1.2 &0.1&144.2  &3 &0.1 &21.1  &2.6  &641.4  &2 &0.2 &47.0  &4.1   &532.6  \\
\hbox{(F2)}   &  &0.1 &10.8&3.6&420.2  &  &0.6 &63.4  &25.0 &948.2  &  &2.3 &2379.9&1127.8&14480.0\\
\hbox{(F3)}   &  &0.0 &1.5 &0.1&177.0  &  &0.1 &22.7  &3.6  &718.2  &  &0.5 &50.7  &4.5   &1077.0 \\
\hbox{(F4)}   &  &0.1 &3.0 &0.1&116.2  &  &0.5 &49.4  &4.8  &548.2  &  &1.9 &172.6 &21.5  &489.8  \\
\hbox{(F5.1)} &  &0.0 &1.1 &0.2&163.4  &  &0.0 &25.8  &4.5  &758.6  &  &0.0 &51.2  &8.4   &597.4  \\
\hbox{(F5.2)} &  &0.0 &1.6 &0.2&169.8  &  &0.0 &31.5  &4.1  &759.8  &  &0.0 &72.5  &16.1  &698.2  \\
\hline
\hbox{(F1)}   &3 &0.0 &2.4 &0.2&349.4  &6 &0.1 &83.0  &15.6 &5006.6 &4 &0.2 &281.5 &11.6  &3955.4 \\
\hbox{(F2)}   &  &0.1 &7.3 &2.3&317.8  &  &0.4 &4.5   &1.7  &113.4  &  &1.9 &265.8 &133.0 &2142.6 \\
\hbox{(F3)}   &  &0.0 &2.5 &0.3&327.8  &  &0.1 &81.9  &10.8 &4357.8 &  &0.2 &263.6 &29.7  &3373.8  \\
\hbox{(F4)}   &  &0.1 &4.5 &0.6&255.0  &  &0.3 &113.1 &54.9 &2986.5 &  &1.6 &813.5 &204.8 &3623.8 \\
\hbox{(F5.1)} &  &0.0 &2.4 &0.4&349.4  &  &0.0 &145.8 &39.9 &5221.0 &  &0.0 &496.7 &117.9 &4619.8 \\
\hbox{(F5.2)} &  &0.0 &3.1 &0.3&373.4  &  &0.0 &20.3  &192.1&5185.4 &  &0.0 &534.9 &130.7 &4085.4 \\
\hline
\hbox{(F1)}   &5 &0.0 &3.9&0.9&847.8   &10&0.0 &177.0 &34.1 &19892.2 &8 &0.2 &2269.1&735.8 &55702.6\\
\hbox{(F2)}   &  &0.1 &2.6&1.1&63.0    &  &0.3 &1.3   &0.4  &35.4    &  &1.3 &16.6  &5.9   &207.8  \\
\hbox{(F3)}   &  &0.0 &4.4&1.3&917.0   &  &0.1 &174.5 &84.1 &17987.0 &  &0.2 &2443.9&357.9 &57277.8\\
\hbox{(F4)}   &  &0.1 &6.0&1.9&599.8   &  &0.2 &233.5 &63.9 &10188.4 &  &1.0 &>1H   &   -  &   -    \\
\hbox{(F5.1)} &  &0.0 &4.9&1.1&878.2   &  &0.0 &340.6 &74.1 &14566.0 &  &0.0 &>1H   &  -   &   -   \\
\hbox{(F5.2)} &  &0.0 &6.9&1.9&869.4   &  &0.0 &610.2 &190.7&21099.0 &  &0.0 &>1H   &    - &   -    \\
\hline
\hbox{(F1)}   &7 &0.0 &4.0&0.7&1153.0  &12&0.0 &219.7 &150.0&30517.8 &10&0.1 &>1H  &   -   &  -  \\
\hbox{(F2)}   &  &0.1 &1.4&0.7&27.4    &  &0.3 &0.7   &0.2  &7.8     &  &1.2 &13.5 & 14.7  &80.6\\
\hbox{(F3)}   &  &0.0 &4.4&0.8&1255.8  &  &0.1 &425.6 &260.4&51375.0 &  &0.2 &>1H  &   -   & - \\
\hbox{(F4)}   &  &0.0 &4.9&1.2&669.8   &  &0.2 &199.5 &147.1&11951.0 &  &0.9 &>1H  &  -    & - \\
\hbox{(F5.1)} &  &0.0 &6.1&2.7&1004.2  &  &0.0 &387.7 &284.1&22344.2 &  &0.0 &>1H  &   -   & - \\
\hbox{(F5.2)} &  &0.0 &9.7&2.6&1399.8  &  &0.0 &964.4 &383.1&37161.4 &  &0.0 &>1H  &  -    & - \\
\hline
\hbox{(F1)}   &10&0.0 &8.5 &1.7&5382.6 &15&0.0 &680.2 &345.5&205131.4 &16&0.1 &3097.7&836.5 &117527.7\\
\hbox{(F2)}   &  &0.0 &0.2 &0.1&3.0    &  &0.2 &1.9   &2.2  &13.0     &  &0.8 &2.0   &0.6   &12.2\\
\hbox{(F3)}   &  &0.0 &9.5 &2.0&5579.8 &  &0.1 &996.8 &25.0 &249648.0 &  &0.2 &>1H   &   -  &  -  \\
\hbox{(F4)}   &  &0.0 &10.5&2.6&3671.0 &  &0.1 &1043.0&18.0 &158639.0 &  &0.6 &>1H   &   -  &  -  \\
\hbox{(F5.1)} &  &0.0 &18.0&5.9&5345.0 &  &0.0 &1701.9&79.0 &208383.0 &  &0.0 &>1H   &   -  &  -  \\
\hbox{(F5.2)} &  &0.0 &26.6&2.6&6025.4 &  &0.0 &>1H  & -  &    -      &  &0.0 &>1H   &   -  &  -  \\
\hline
\hbox{(F1)}   &12&0.0 &16.5&1.8 &16316.2 &22&0.0 &2453.0&5.0  & 2\cdot 10^6 &20&0.1 &>1H  &   -  &   -  \\
\hbox{(F2)}   &  &0.0 &0.1 &0.0 &1.0     &  &0.0 &0.2   &0.1  &1.0       &  &0.6 &3.8  &6.6   &1.8\\
\hbox{(F3)}   &  &0.0 &16.4&1.5 &14263.0 &  &0.1 &3298.1&150.8& 2\cdot 10^6 &  &0.2 &>1H  &   -  &  -   \\
\hbox{(F4)}   &  &0.0 &21.6&4.0 &13839.0 &  &0.1 &3107.5&5.4  & 2\cdot 10^6 &  &0.4 &>1H  &   -  &  -   \\
\hbox{(F5.1)} &  &0.0 &33.2&11.9&15626.6 &  &0.0 &>1H  &   -  &     -    &  &0.0 &>1H  &   -  &  -   \\
\hbox{(F5.2)} &  &0.0 &46.5&6.7 &14411.4 &  &0.0 &>1H  &   -  &     -    &  &0.0 &>1H  &   -  &  -    \\
\hline
\end{array}$$
\caption{\label{freeSS} Customers prefer closer sites and
self-service is allowed}
\end{table}

\begin{table}
$$
\hspace{-2cm}
\begin{array}{|c|crrrr|rrrrr|rrrrr|}
\hline \multicolumn{1}{|c}{} & \multicolumn{5}{|c|}{M=20} &
\multicolumn{5}{|c|}{M=30} &
\multicolumn{5}{|c|}{M=40} \\
\hline
  & p & \multicolumn{1}{c}{LP} & \multicolumn{1}{c}{\bar{t}} & \multicolumn{1}{c}{\sigma_t} &
\multicolumn{1}{c}{n} & \multicolumn{1}{|c}{p} & \multicolumn{1}{c}{LP} & \multicolumn{1}{c}{\bar{t}} & \multicolumn{1}{c}{\sigma_t} & \multicolumn{1}{c}{n}  & \multicolumn{1}{|c}{p} & \multicolumn{1}{c}{LP} & \multicolumn{1}{c}{\bar{t}} & \multicolumn{1}{c}{\sigma_t} & \multicolumn{1}{c|}{n} \\
\hline
\hbox{(F1)}   &2 &0.0 &1.1 &0.2 &120.6  &3 &0.1&23.9 &2.7  &799.0     &2 &0.2 &44.9  &9.8   &437.0           \\
\hbox{(F2)}   &  &0.1 &28.2&19.0&2117.8 &  &0.5&814.2&293.8&21899.4   &  &2.0 &2717.5&855.7 &17868.5       \\
\hbox{(F3)}   &  &0.0 &1.4 &0.3 &167.8  &  &0.1&26.5 &2.8  &861.9     &  &0.3 &53.1  &11.2  &482.2          \\
\hbox{(F4)}   &  &0.1 &3.6 &0.9 &126.6  &  &0.4&55.1 &5.4  &644.2     &  &1.8 &202.9 &59.6  &518.2         \\
\hbox{(F5.1)} &  &0.0 &1.0 &0.4 &127.8  &  &0.0&28.1 &4.4  &772.2     &  &0.0 &48.3  &9.5   &469.8         \\
\hbox{(F5.2)} &  &0.0 &1.7 &0.3 &187.0  &  &0.0&39.7 &6.9  &977.8     &  &0.0 &71.0  &16.1  &575.0           \\
\hline
\hbox{(F1)}   &3 &0.0 &0.4 &0.2 &129.4  &6 &0.1&52.1 &12.2 &2672.6    &4 &0.2 &312.0 &99.7  &4173.4         \\
\hbox{(F2)}   &  &0.1 &16.6&7.9 &1791.0 &  &0.4&52.0 &26.3 &2749.8    &  &1.8 &>1H   &   -  &   -           \\
\hbox{(F3)}   &  &0.0 &2.5 &0.4 &315.4  &  &0.1&66.9 &25.0 &3258.6    &  &0.2 &326.3 &84.6  &4029.4         \\
\hbox{(F4)}   &  &0.1 &4.5 &0.4 &213.4  &  &0.3&109.4&46.9 &1859.0    &  &1.4 &922.7 &179.8 &3653.8         \\
\hbox{(F5.1)} &  &0.0 &2.4 &0.4 &329.8  &  &0.0&92.0 &51.0 &2789.0    &  &0.0 &521.6 &92.8  &4542.2         \\
\hbox{(F5.2)} &  &0.0 &3.1 &0.4 &324.2  &  &0.0&141.2&50.0 &3219.8    &  &0.0 &656.4 &168.2 &4777.8         \\
\hline
\hbox{(F1)}   &5 &0.0 &2.9 &0.7 &584.6  &10&0.1&87.7 &37.4 &7987.8    &8 &0.2 &2000.2&657.0 &43647.4        \\
\hbox{(F2)}   &  &0.1 &2.5 &1.0 &364.6  &  &0.4&4.5  &2.3  &420.4     &  &1.4 &846.7 &493.9 &28928.2          \\
\hbox{(F3)}   &  &0.1 &3.2 &0.5 &597.0  &  &0.1&82.3 &6.6  &6807.0    &  &0.2 &2105.2&302.5 &38412.2       \\
\hbox{(F4)}   &  &0.0 &5.5 &1.2 &453.4  &  &0.2&116.3&18.9 &4276.0    &  &1.2 &3118.8&411.1 &21404.5       \\
\hbox{(F5.1)} &  &0.0 &4.3 &1.2 &690.6  &  &0.0&231.7&16.2 &9669.0    &  &0.0 &>1H  &   -   &   -           \\
\hbox{(F5.2)} &  &0.0 &6.0 &1.6 &751.0  &  &0.0&298.3&16.0 &9997.8    &  &0.0 &>1H  &   -   &   -           \\
\hline
\hbox{(F1)}   &7 &0.0 &3.0 &0.7 &722.2  &12&0.0&49.6 &23.9 &5036.2    &10&0.2 &2624.0&694.7 &62600.0\\
\hbox{(F2)}   &  &0.1 &1.0 &0.6 &101.4  &  &0.3&1.7  &0.4  &97.4      &  &1.3 &139.6 &112.2 &4933.4\\
\hbox{(F3)}   &  &0.0 &2.5 &1.1 &595.4  &  &0.1&48.8 &17.7 &4261.0    &  &0.2 &2763.5&449.3 &67043.0\\
\hbox{(F4)}   &  &0.0 &4.0 &2.2 &489.8  &  &0.2&50.8 &13.8 &1902.6    &  &0.9 &>1H   &  -    &  -   \\
\hbox{(F5.1)} &  &0.0 &4.7 &2.3 &908.6  &  &0.0&115.0&56.5 &5427.4    &  &0.0 &>1H   &  -    &  -   \\
\hbox{(F5.2)} &  &0.0 &5.8 &3.0 &918.6  &  &0.0&215.6&124.4&6569.0    &  &0.0 &>1H   &  -    &  -   \\
\hline
\hbox{(F1)}   &10&0.0 &1.0 &0.5 &315.4  &15&0.0&16.1 &5.8  &2396.6    &16&0.1 &1555.1&1205.4&41557.0       \\
\hbox{(F2)}   &  &0.1 &0.4 &0.2 &28.6   &  &0.2&2.2  &2.0  &21.8      &  &0.8 &4.2   &0.9   &136.0         \\
\hbox{(F3)}   &  &0.0 &1.3 &0.7 &411.0  &  &0.1&15.8 &0.7  &2054.0    &  &0.2 &1677.0&1147.9&43331.5       \\
\hbox{(F4)}   &  &0.0 &1.5 &0.8 &207.4  &  &0.2&13.8 &4.7  &647.0     &  &0.6 &1685.2&1138.8&17580.5       \\
\hbox{(F5.1)} &  &0.0 &1.5 &1.2 &350.6  &  &0.0&27.7 &1.3  &1697.0    &  &0.0 &>1H   &   -   &   -          \\
\hbox{(F5.2)} &  &0.0 &2.0 &1.2 &374.6  &  &0.0&29.7 &1.3  &1710.2    &  &0.0 &>1H   &   -   &   -          \\
\hline
\hbox{(F1)}   &12&0.0 &0.4 &0.2 &129.4  &22&0.0&0.1  &0.0 &2.6        &20&0.1 &315.3 &302.2 &19865.8     \\
\hbox{(F2)}   &  &0.0 &0.4 &0.1 &10.6   &  &0.1&0.2  &0.2 &1.8        &  &0.6 &2.0   &0.4   &53.8         \\
\hbox{(F3)}   &  &0.0 &0.4 &0.2 &107.4  &  &0.0&0.2  &0.0 &3.0        &  &0.2 &127.5 &26.7  &6387.0           \\
\hbox{(F4)}   &  &0.0 &0.6 &0.3 &87.4   &  &0.1&0.3  &0.0 &3.0        &  &0.4 &137.6 &61.7  &905.10             \\
\hbox{(F5.1)} &  &0.0 &0.4 &0.2 &114.6  &  &0.0&0.1  &0.0 &3.4        &  &0.0 &353.3 &364.8 &7442.0             \\
\hbox{(F5.2)} &  &0.0 &0.6 &0.3 &124.6  &  &0.0&0.1  &0.0 &3.0        &  &0.0 &412.4 &371.2 &7987.2             \\
\hline
\end{array}
$$
\caption{\label{nofreeSS} Customers prefer closer sites but
self-service is forbidden}
\end{table}

\begin{table}
$$
\hspace{-2cm}
\begin{array}{|c|crrrr|rrrrr|rrrrr|}
\hline \multicolumn{1}{|c}{} & \multicolumn{5}{|c|}{M=20} &
\multicolumn{5}{|c|}{M=30} &
\multicolumn{5}{|c|}{M=40} \\
\hline
  & p & \multicolumn{1}{c}{LP} & \multicolumn{1}{c}{\bar{t}} & \multicolumn{1}{c}{\sigma_t} &
\multicolumn{1}{c}{n} & \multicolumn{1}{|c}{p} & \multicolumn{1}{c}{LP} & \multicolumn{1}{c}{\bar{t}} & \multicolumn{1}{c}{\sigma_t} & \multicolumn{1}{c}{n}  & \multicolumn{1}{|c}{p} & \multicolumn{1}{c}{LP} & \multicolumn{1}{c}{\bar{t}} & \multicolumn{1}{c}{\sigma_t} & \multicolumn{1}{c|}{n} \\
\hline
\hbox{(F1)}   &2 &0.0&1.6 &0.1  &201.8      &3 &0.1&28.1 &4.1 &869.8       &2 &0.3&90.2  &8.6  &1082.6\\
\hbox{(F2)}   &  &0.1&62.7&26.7 &4650.6     &  &0.6&>1H  &  - &   -        &  &2.1&>1H   &   - &  -   \\
\hbox{(F3)}   &  &0.0&2.0 &0.3  &242.6      &  &0.1&34.0 &8.0 &1033.4      &  &0.0&96.6  &13.0 &1150.2 \\
\hbox{(F4)}   &  &0.1&3.9 &0.6  &153.8      &  &0.4&72.5 &7.8 &756.2       &  &0.0&342.2 &12.2 &1095.0  \\
\hbox{(F5.1)} &  &0.0&1.5 &0.4  &191.8      &  &0.0&42.7 &4.1 &1155.4      &  &0.0&100.3 &14.5 &1095.0 \\
\hbox{(F5.2)} &  &0.0&1.7 &0.5  &205.0      &  &0.0&47.1 &6.8 &1000.6      &  &0.0&92.5  &6.6  &1183.8\\
\hline
\hbox{(F1)}   &3 &0.0&2.2 &0.1  &273.0      &6 &0.1&66.4 &17.7&2903.0      &4 &0.2&459.6 &85.5  &5656.6\\
\hbox{(F2)}   &  &0.1&15.5&4.3  &1262.6     &  &0.4&113.5&75.0&5463.8      &  &2.0&>1H   &   -  &  -  \\
\hbox{(F3)}   &  &0.0&2.5 &0.1  &293.4      &  &0.1&63.2 &10.6&2436.6      &  &0.3&507.0 &57.1  &5660.2\\
\hbox{(F4)}   &  &0.1&5.1 &0.2  &243.8      &  &0.3&125.3&15.9&1989.4      &  &1.5&1236.3&254.0 &4778.6\\
\hbox{(F5.1)} &  &0.0&74.8&161.1&264.1      &  &0.0&118.2&27.5&3243.8      &  &0.0&917.3 &114.8 &6964.2\\
\hbox{(F5.2)} &  &0.0&2.9 &0.4  &283.0      &  &0.0&125.6&34.8&2409.4      &  &0.0&859.3 &193.8 &5883.8\\
\hline
\hbox{(F1)}   &5 &0.0&2.5 &0.2  &389.8      &10&0.1&32.0&13.5  &2085.8      &8 &0.2&977.1 &633.9 &15782.6\\
\hbox{(F2)}   &  &0.1&2.4 &0.7  &329.4      &  &0.3&2.7 &0.4   &137.8       &  &1.3&212.8 &234.4 &4813.8\\
\hbox{(F3)}   &  &0.0&2.3 &0.4  &346.2      &  &0.1&27.8&0.1   &1443.0      &  &0.2&862.9 &402.6 &12009.4\\
\hbox{(F4)}   &  &0.1&4.3 &0.9  &299.4      &  &0.2&58.5&13.4  &1383.0      &  &0.9&1496.8&762.3 &7941.4\\
\hbox{(F5.1)} &  &0.0&2.8 &0.5  &394.2      &  &0.0&43.4&29.2  &1684.0      &  &0.0&1969.2&879.5 &15085.4\\
\hbox{(F5.2)} &  &0.0&4.0 &1.1  &456.2      &  &0.0&45.6&29.1  &1710.2      &  &0.0&2169.6&825.9 &12589.8\\
\hline
\hbox{(F1)}   &7 &0.0&1.4 &0.4  &270.2      &12&0.0&18.6 &15.6&1384.6      &10&0.2&985.2 &456.9 &19158.2\\
\hbox{(F2)}   &  &0.1&1.1 &0.8  &44.6       &  &0.3&1.1  &0.2 &39.4        &  &1.2&28.4  &13.1  &645.4\\
\hbox{(F3)}   &  &0.0&1.7 &0.6  &329.4      &  &0.1&20.5 &16.3&1495.4      &  &0.2&651.4 &139.3 &11759.0\\
\hbox{(F4)}   &  &0.0&2.0 &0.5  &175.8      &  &0.2&23.8 &21.1&721.4       &  &0.9&1234.7&145.2 &7846.0\\
\hbox{(F5.1)} &  &0.0&1.8 &0.5  &326.2      &  &0.0&36.7  &38.1&1661.4     &  &0.0&2390.3&521.4 &21169.0\\
\hbox{(F5.2)} &  &0.0&2.2 &1.2  &286.2      &  &0.0&57.5 &27.1&1727.4      &  &0.0&2675.9&546.3 &22998.2\\
\hline
\hbox{(F1)}   &10&0.0&0.5 &0.3  &145.0      &15&0.0&4.7  &4.7 &50.0        &16&0.1&280.2&202.5  &9982.2\\
\hbox{(F2)}   &  &0.0&0.5 &0.3  &16.2       &  &0.2&0.8  &0.1 &19.0        &  &0.8&2.7  &0.5    &31.8\\
\hbox{(F3)}   &  &0.0&0.6 &0.4  &166.2      &  &0.1&6.4  &7.1 &648.0       &  &0.2&336.1&202.0  &12150.0\\
\hbox{(F4)}   &  &0.0&0.7 &0.2  &70.2       &  &0.1&6.9  &6.3 &321.0       &  &0.7&315.8&118.8  &4096.0\\
\hbox{(F5.1)} &  &0.0&0.6 &0.3  &155.8      &  &0.0&6.2  &8.0 &451.0       &  &0.0&798.8&256.2  &13713.0\\
\hbox{(F5.2)} &  &0.0&0.9 &0.7  &185.4      &  &0.0&7.2  &8.1 &461.5       &  &0.0&899.4&265.2  &13942.0\\
\hline
\hbox{(F1)}   &12&0.0&0.2 &0.1  &33.0       &22&0.0&0.0  &0.1 &2.0         &20&0.0&75.3 &63.6    &3708.2\\
\hbox{(F2)}   &  &0.0&0.3 &0.1  &6.2        &  &0.1&0.1  &0.0 &1.0         &  &0.6&2.1  &0.2     &39.4 \\
\hbox{(F3)}   &  &0.0&0.2 &0.2  &51.0       &  &0.0&0.2  &0.1 &3.0         &  &0.2&71.1 &8.3     &3689.0 \\
\hbox{(F4)}   &  &0.0&0.4 &0.2  &27.0       &  &0.1&0.3  &0.1 &3.0         &  &0.5&96.4 &40.9    &1975.0\\
\hbox{(F5.1)} &  &0.0&0.2 &0.1  &51.8       &  &0.0&0.1  &0.0 &3.0         &  &0.0&108.8&36.3    &2800.0\\
\hbox{(F5.2)} &  &0.0&0.3 &0.2  &49.0       &  &0.0&0.1  &0.0 &3.0         &  &0.0&158.3&37.2    &2901.2\\
\hline
\end{array}
$$
\caption{\label{random} Random preferences}
\end{table}

\section{Improving formulations}
\label{sec:improving}

Some of these formulations have been improved and solution method
have been again computationally compared on the same testbed.

Table \ref{f1mejor} reports a comparative analysis of the results
provided when solving formulation (F1) with and without a separation
method.

\begin{table}
$$
\hspace{-2.5cm}
\begin{array}{|c|rrrrr|rrrrr|rrrrr|}
\hline
 &\multicolumn{5}{|c|}{M=20} & \multicolumn{5}{|c|}{M=30} & \multicolumn{5}{|c|}{M=40} \\
\hline
 & \multicolumn{1}{c}{p} & \multicolumn{1}{c}{LP} & \multicolumn{1}{c}{\bar{t}} &
 \multicolumn{1}{c}{\sigma_t } &
 \multicolumn{1}{c|}{n} & \multicolumn{1}{c}{p} & \multicolumn{1}{c}{LP} & \multicolumn{1}{c}{\bar{t}} &
\multicolumn{1}{c}{\sigma_t} & \multicolumn{1}{c|}{n} &
\multicolumn{1}{c}{p} & \multicolumn{1}{c}{LP} &
\multicolumn{1}{c}{\bar{t}} & \multicolumn{1}{c}{\sigma_t} & \multicolumn{1}{c|}{n} \\
\hline
  \multicolumn{16}{|c|}{\hbox{Customers prefer closer sites and
self-service is allowed}} \\
\hline
\hbox{(F1)}   &2 &0.0 &1.2 &0.1 &144.2  & 3 &0.1 &21.1  &2.6  &641.4    &    2 &0.2 &47.0  &4.1   &532.6   \\
\hbox{(F1R)}  &  &0.0 &1.7 &0.1 &155.0  &   &0.1 &24.8  &2.8  &557.0    &      &0.2 &74.4  &10.8  &496.6   \\
\hline
\hbox{(F1)}   &3 &0.0 &2.4 &0.2 &349.4  & 6 &0.1 &83.0  &15.6 &5006.6   &    4 &0.2 &281.5 &11.6  &3955.4  \\
\hbox{(F1R)}  &  &0.0 &2.8 &0.2 &285.0  &   &0.1 &95.7  &23.9 &3910.2   &      &0.2 &359.1 &65.2  &3355.8  \\
\hline
\hbox{(F1)}   &5 &0.0 &3.9 &0.9 &847.8  & 10&0.0 &177.0 &34.1 &19892.2  &    8 &0.2 &2269.1&735.8 &55702.6 \\
\hbox{(F1R)}  &  &0.0 &4.5 &1.4 &622.6  &   &0.1 &180.1 &56.3 &11339.8  &      &0.1 &2655.5&645.9 &42317.8 \\
\hline
\hbox{(F1)}   &7 &0.0 &4.0 &0.7 &1153.0 & 12&0.0 &219.7 &150.0&30517.8  &    10&0.1 &>1H  &   - &   -     \\
\hbox{(F1R)}  &  &0.0 &4.9 &1.1 &786.2  &   &0.0 &188.5 &113.1&13592.6  &      &0.1 &>1H  &   - &   -     \\
\hline
\hbox{(F1)}   &10&0.0 &8.5 &1.7 &5382.6 & 15&0.0 &680.2 &345.5&205131.4 &    16&0.1 &3097.7 &836.5 &117527.7\\
\hbox{(F1R)}  &  &0.0 &12.9&2.4 &4125.8 &   &0.0 &635.6 &89.0 &116827.0 &      &0.1 &>1H  &  -  &   -   \\
\hline
\hbox{(F1)}   &12&0.0 &16.5&1.8 &16316.2& 22&0.0 &2453.0&5.0  &2\cdot 10^6&    20&0.1 &>1H  &  -  &   -  \\
\hbox{(F1R)}  &  &0.0 &25.3&3.3 &14071.8&   &0.0 &>1H & -   &    -    &      &0.1 &>1H  &  -  &   -   \\
\hline
  \multicolumn{16}{|c|}{\hbox{Customers prefer closer sites but
self-service is forbidden}} \\
\hline
\hbox{(F1)}   &2 &0.0&1.1  &0.2  &120.6 & 3 &0.1 &23.9 &2.7  &799.0 &  2 &0.2 &44.9  &9.8   &437.0   \\
\hbox{(F1R)}  &  &0.0&1.6  &0.2  &120.6 &   &0.1 &30.4 &4.7  &655.8 &    &0.2 &64.9  &4.2   &353.8 \\
\hline
\hbox{(F1)}   &3 &0.0&0.4  &0.2  &129.4 & 6 &0.1 &52.1 &12.1 &2672.6&  4 &0.2 &312.0 &99.7  &4173.4 \\
\hbox{(F1R)}  &  &0.0&2.7  &0.4  &245.8 &   &0.1 &75.6 &24.6 &2389.0&    &0.2 &416.1 &88.6  &3450.2\\
\hline
\hbox{(F1)}   &5 &0.0&2.9  &0.7  &584.6 & 10&0.1 &87.7 &37.4 &7987.8&  8 &0.2 &2000.2&657.0 &43647.4 \\
\hbox{(F1R)}  &  &0.0&4.1  &0.9  &499.8 &   &0.1 &87.5 &36.2 &4195.0&    &0.2 &2405.5&394.5 &29669.0\\
\hline
\hbox{(F1)}   &7 &0.0&3.0  &0.7  &722.2 & 12&0.0 &49.6 &23.9 &5036.2&  10&0.2 &2624.0&694.7 &62600.0\\
\hbox{(F1R)}  &  &0.0&3.8  &1.1  &519.0 &   &0.1 &61.4 &29.8 &2890.2&    &0.1 &3020.4&675.4 &35502.0\\
\hline
\hbox{(F1)}   &10&0.0&1.0  &0.5  &315.4 & 15&0.0 &16.1 &5.8  &2396.6&  16&0.1 &1555.1&1205.4&41557.0 \\
\hbox{(F1R)}  &  &0.0&2.1  &1.1  &337.8 &   &0.1 &23.4 &7.7  &1117.8&    &0.2 &1261.4&665.5 &24426.3\\
\hline
\hbox{(F1)}   &12&0.0&0.4  &0.2  &129.4 & 22&0.0 &0.1  &0.0  &2.6   &  20&0.1 &315.3 &302.2 &19865.8 \\
\hbox{(F1R)}  &  &0.0&0.7  &0.4  &90.6  &   &0.0 &0.7  &0.4  &2.6   &    &0.1 &523.3 &286.8 &11281.8 \\
\hline
 \multicolumn{16}{|c|}{\hbox{Random preferences}} \\
\hline
\hbox{(F1)}   &2 &0.0&1.6 &0.1  &201.8 & 3 &0.1 &28.1 &4.1  &869.8 &   2 &0.3 &90.2  &8.6   &1082.6\\
\hbox{(F1R)}  &  &0.0&2.1 &0.6  &183.4 &   &0.1 &40.4 &1.6  &728.6 &     &0.2 &147.9 &13.0  &1052.3\\
\hline
\hbox{(F1)}   &3 &0.0&2.2 &0.1  &273.0 & 6 &0.1 &66.4 &17.7 &2903.0&   4 &0.2 &459.6 &85.5  &5656.6\\
\hbox{(F1R)}  &  &0.0&2.8 &0.4  &218.2 &   &0.1 &84.4 &12.1 &2320.6&     &0.2 &767.3 &93.0  &5431.0\\
\hline
\hbox{(F1)}   &5 &0.0&2.5 &0.2  &389.8 & 10&0.1 &32.0 &13.5 &2085.8&   8 &0.2 &977.1 &633.9 &15782.6\\
\hbox{(F1R)}  &  &0.0&3.5 &0.9  &381.0 &   &0.1 &43.2 &8.3  &1470.2&     &0.3 &1033.7&464.9 &9091.4\\
\hline
\hbox{(F1)}   &7 &0.0&1.4 &0.4  &270.2 & 12&0.0 &18.6 &15.6 &1384.6&   10&0.2 &985.2 &456.9 &19158.2\\
\hbox{(F1R)}  &  &0.0&2.2 &0.7  &248.6 &   &0.1 &30.4 &21.8 &1037.0&     &0.2 &1118.1&1118.7&14588.3\\
\hline
\hbox{(F1)}   &10&0.0&0.5 &0.3  &145.0 & 15&0.0 &4.7  &4.7  &50.0  &   16&0.1 &280.2 &202.5 &9982.2\\
\hbox{(F1R)}  &  &0.0&1.0 &0.5  &128.2 &   &0.1 &8.7  &9.4  &372.6 &     &0.2 &390.0 &209.2 &4941.0\\
\hline
\hbox{(F1)}   &12&0.0&0.2 &0.1  &33.0  & 22&0.0 &0.0  &0.1  &2.0   &   20&0.0 &75.3  &63.6  &3708.2  \\
\hbox{(F1R)}  &  &0.0&0.3 &0.2  &32.2  &   &0.0 &0.2  &0.0  &3.8   &     &0.2 &124.6 &47.2  &2098.2\\
\hline
\end{array}$$
\caption{\label{f1mejor} Comparison of (F1) with and without
additional constraints}
\end{table}

Table \ref{VI} shows the running times when solving formulation (F2)
with several sets of additional inequalities. A separation procedure
have been implemented and the results are given in Table \ref{KK}.
The results of solving the preprocessed formulation (F2) are
presented in Table \ref{f2Pre}.
\begin{eqnarray}
\label{p101} & & x_{ik}\ge x_{i,k+1}, \quad \forall i\in A, k=2,\ldots,M-p;  \\
\label{p102} & & pz_{ik} \le \sum_{j:\ O_{ij}\ge k } y_j, \quad \forall i\in A, k=2,\ldots,M-p; \\
\label{p103} & & (p-1)(z_{ik}-z_{i,k+1}) \le \sum_{j:\ O_{ij}\ge k+1} y_j, \quad \forall i\in A, k=2,\ldots,M-p; \\
\label{p104} & & \sum_{i=1}^s x_{M+1-i,k}\ge \sum_{i\in S} z_{ik}, \quad \forall s\in A,S\subset A:\ |S|=s,  k=2,\ldots,M-p+1; \\
\label{p105} & & z_{ik}+\sum_{j:\ O_{ij}<k} y_j \ge 1, \quad \forall
i\in A, k=2,\ldots,M-p+1.
\end{eqnarray}

\begin{table}
{\small
$$
\hspace{-2cm}
\begin{array}{|c|rr|rr|rr|rr|rr|rr|}
\hline \multicolumn{1}{|c}{\hbox{Ineq.}} &
\multicolumn{1}{|c}{\bar{t}} & \multicolumn{1}{c|}{n} &
\multicolumn{1}{|c}{\bar{t}} & \multicolumn{1}{c|}{n} &
\multicolumn{1}{|c}{\bar{t}} & \multicolumn{1}{c|}{n} &
\multicolumn{1}{|c}{\bar{t}} & \multicolumn{1}{c|}{n} &
\multicolumn{1}{|c}{\bar{t}} & \multicolumn{1}{c|}{n} &
\multicolumn{1}{|c}{\bar{t}} &
\multicolumn{1}{c|}{n} \\
\hline
 \multicolumn{13}{|c|}{M=20} \\
\hline \multicolumn{1}{|c}{} & \multicolumn{2}{|c|}{p=2} &
\multicolumn{2}{|c|}{p=3} &
\multicolumn{2}{|c|}{p=5}        & \multicolumn{2}{|c|}{p=7} & \multicolumn{2}{|c|}{p=10} & \multicolumn{2}{|c|}{p=12} \\
\hline \hbox{None}
 &10.8&420.2  &7.3 &317.8&   2.6&63.0&   1.4&27.4&   0.2&3.0     &0.1&1.0    \\
(\ref{p101})
 &4.6 &96.2   &2.9 &71.0 &   1.2&42.2&   0.4&10.0&   0.2&5.0     &0.1&1.0    \\
(\ref{p102})
 &13.8&436.2  &8.6 &261.8&   3.1&57.4&   2.2&25.4&   0.4&5.0     &0.2&5.0    \\
(\ref{p103})
 &14.1&415.0  &9.5 &370.6&   3.1&68.6&   1.7&15.4&   0.1&1.0     &0.2&4.2    \\
(\ref{p105})
 &11.5&386.7  &8.5 &284.2&   3.2&49.8&   1.5&18.2&   0.3&5.4     &0.1&1.0    \\
(\ref{p101})(\ref{p102})
 &4.7 &89.8   &2.6 &81.8 &   1.5&31.0&   0.7&9.4&   0.3&5.0     &0.1&1.0    \\
(\ref{p101})(\ref{p103})
 &5.5 &71.4   &3.1 &91.8 &   1.6&29.4&   0.9&8.2 &   0.6&10.6    &0.1&1.0    \\
(\ref{p101})(\ref{p105})
 &5.1 &101.0  &2.9 &69.8 &   1.1&23.4&   1.1&14.2&   0.4&8.6     &0.1&1.0    \\
(\ref{p102})(\ref{p103})
 &15.4&453.4  &7.1 &173.8&   2.2&32.6&   1.8&17.0&   0.5&9.0     &0.1&1.0    \\
(\ref{p102})(\ref{p105})
 &13.4&294.2  &8.3 &221.0&   2.4&35.0&   2.0&12.6&   0.3&1.0     &0.1&1.0    \\
(\ref{p103})(\ref{p105})
 &14.5&371.8  &8.2 &199.8&   2.8&38.2&   1.6&26.6&   0.2&3.4     &0.1&1.0    \\
(\ref{p101})(\ref{p102})(\ref{p103})
 &5.0 &111.4  &3.2 &75.8 &   1.0&23.4&   0.9&11.8&   0.1&1.0     &0.1&1.0    \\
(\ref{p101})(\ref{p102})(\ref{p105})
 &6.1 &99.0   &3.6 &75.0 &   0.9&34.2&   0.7&6.2 &   0.3&4.2     &0.1&1.0    \\
(\ref{p101})(\ref{p103})(\ref{p105})
 &7.2 &114.6  &5.2 &98.2 &   1.0&37.0&   1.3&9.4 &   0.2&4.6     &0.2&4.2    \\
(\ref{p102})(\ref{p103})(\ref{p105})
 &18.2&402.6  &13.3&297.0&   3.5&52.6&   2.5&26.6&   0.9&13.0    &0.1&1.0    \\
(\ref{p101})(\ref{p102})(\ref{p103})(\ref{p105})
 &7.0 &97.4   &5.1 &85.4 &   1.6&38.6&   1.0&9.0 &   0.5&8.2     &0.3&4.2    \\
\hline
 \multicolumn{13}{|c|}{M=30} \\
\hline \multicolumn{1}{|c}{} & \multicolumn{2}{|c|}{p=3} &
\multicolumn{2}{|c|}{p=6} &
\multicolumn{2}{|c|}{p=10}       & \multicolumn{2}{|c|}{p=12} & \multicolumn{2}{|c|}{p=15} & \multicolumn{2}{|c|}{p=22} \\
\hline \hbox{None}
  &63.4&948.2   &4.5 &113.4&   1.3 &35.4 &   0.7&7.8    &1.9&13.0      & 0.2&1.0  \\
(\ref{p101})
  &10.0&159.0   &3.0 &8.7 &   1.1 &15.0&   0.7&2.3   &1.4&11.4     & 0.2&1.0  \\
(\ref{p102})
  &55.8&630.6   &8.2 &24.1&   1.4 &25.4&   0.9&10.8  &0.3&1.0      & 0.2&1.0  \\
(\ref{p103})
  &80.8&1005.0  &5.8 &17.6&   1.8 &33.4&   0.7&1.8   &1.1&7.0      & 0.2&1.0  \\
(\ref{p105})
  &68.9&984.6   &5.4 &52.2&   1.3 &25.0&   0.7&4.9   &1.2&7.0      & 0.2&1.0  \\
(\ref{p101})(\ref{p102})
  &11.7&117.8   &9.3 &24.8&   1.9 &22.2&   1.1&3.0   &0.6&1.0      & 0.2&1.0  \\
(\ref{p101})(\ref{p103})
  &14.3&178.6   &4.1 &26.9&   1.6 &27.0&   2.2&3.7   &1.2&6.6      & 0.2&1.0  \\
(\ref{p101})(\ref{p105})
  &13.5&143.8   &5.2 &58.6&   2.1 &16.2&   1.1&3.0   &1.0&5.8      & 0.2&1.0  \\
(\ref{p102})(\ref{p103})
  &77.1&671.8   &5.7 &23.9&   1.6 &21.4&   0.7&4.1   &0.4&1.0      & 0.3&1.0  \\
(\ref{p102})(\ref{p105})
  &75.7&675.0   &8.7 &48.3&   2.3 &30.6&   1.2&1.7   &0.3&1.0      & 0.3&1.0  \\
(\ref{p103})(\ref{p105})
  &91.7&852.6   &7.4 &86.8&   1.5 &17.0&   0.8&3.3   &1.3&7.0      & 0.2&1.0  \\
(\ref{p101})(\ref{p102})(\ref{p103})
  &15.5&141.0   &7.2 &23.4&   1.7 &19.4&   2.0&2.4   &0.5&1.0      & 0.3&1.0  \\
(\ref{p101})(\ref{p102})(\ref{p105})
  &21.4&176.6   &4.8 &22.8&   2.5 &35.8&   1.4&2.4   &0.4&1.0      & 0.3&1.0  \\
(\ref{p101})(\ref{p103})(\ref{p105})
  &15.6&105.4   &4.7 &38.1&   1.7 &15.8&   1.2&4.3   &0.8&5.0      & 0.3&1.0  \\
(\ref{p102})(\ref{p103})(\ref{p105})
  &97.9&777.4   &12.0&34.4&   2.1 &22.6&   1.2&1.7   &0.4&1.0      & 0.2&1.0  \\
(\ref{p101})(\ref{p102})(\ref{p103})(\ref{p105})
  &23.4&169.8   &16.7&28.5&   2.7 &20.2&   1.8&2.4   &1.0&5.0      & 0.3&1.0  \\
\hline
 \multicolumn{13}{|c|}{M=40} \\
\hline \multicolumn{1}{|c}{} & \multicolumn{2}{|c|}{p=2} &
\multicolumn{2}{|c|}{p=4} & \multicolumn{2}{|c|}{p=8}        &
\multicolumn{2}{|c|}{p=10} & \multicolumn{2}{|c|}{p=16} &
\multicolumn{2}{|c|}{p=20} \\
\hline \hbox{None}
  & 2379.9&14480.0     &265.8 &2142.6 &   16.6&207.8 &   13.5&80.6  &   2.0 &12.2   &3.8&1.8\\
(\ref{p101})
  & 176.9 &838.2       &51.6  &413.8 &   13.1&143.8&   7.2 &86.5 &   2.3 &18.5   &1.0&1.0 \\
(\ref{p102})
  & 2106.9&7885.4      &283.1 &1449.0&   32.6&257.8&   9.0 &93.5 &   2.0 &7.0    &0.6&1.0 \\
(\ref{p103})
  & >1H  &   -         &359.3 &2138.2&   31.4&189.8&   9.8 &92.0 &   2.0 &8.5    &4.3&8.0 \\
(\ref{p105})
  & >1H  &    -     &338.4 &1865.8&   18.8&159.4&   10.6&108.0&   1.7 &7.5    &0.6&1.0 \\
(\ref{p101})(\ref{p102})
  & 185.7 &795.8       &66.0  &348.6 &   21.0&98.2 &   9.4 &93.5 &   2.4 &16.5   &1.0&1.0 \\
(\ref{p101})(\ref{p103})
  & 165.4 &562.2       &58.8  &374.2 &   25.6&110.2&   7.0 &66.5 &   2.1 &4.0    &1.0&1.0 \\
(\ref{p101})(\ref{p105})
  & 319.1 &642.6       &66.4  &423.8 &   26.0&120.2&   7.0 &73.0 &   2.2 &8.5    &0.6&1.0 \\
(\ref{p102})(\ref{p103})
  & 2251.6&6124.5      &766.1 &4007.0&   33.8&182.2&   30.3&106.5&   8.9 &7.5    &1.0&1.0 \\
(\ref{p102})(\ref{p105})
  & 2487.3&7410.5      &356.1 &1290.2&   32.9&232.6&   19.1&102.5&   10.7&7.0    &1.0&1.0 \\
(\ref{p103})(\ref{p105})
  & >1H  &   -     &1022.1&1708.0&   48.7&208.2&   13.1&99.5 &   8.9 &9.0    &1.6&2.0 \\
(\ref{p101})(\ref{p102})(\ref{p103})
  & 219.4 &532.6       &78.4  &333.8 &   46.6&155.0&   18.5&72.5 &   9.2 &13.5   &1.5&1.0 \\
(\ref{p101})(\ref{p102})(\ref{p105})
  & 230.2 &579.4       &115.7 &509.8 &   34.3&87.0 &   23.2&85.0 &   6.3 &4.5    &1.6&1.0 \\
(\ref{p101})(\ref{p103})(\ref{p105})
  & 320.0 &790.2       &122.5 &399.0 &   47.5&117.4&   14.1&107.0&   5.5 &4.5    &2.1&1.0 \\
(\ref{p102})(\ref{p103})(\ref{p105})
  & 2700.5&5293.0      &667.8 &1988.6&   65.5&158.6&   23.2&85.0 &   14.6&8.5    &5.2&1.0 \\
(\ref{p101})(\ref{p102})(\ref{p103})(\ref{p105})
  &353.3  &718.0       &155.2 &473.4 &   73.6&119.0&   29.7&72.0 &   8.6 &4.0    &5.8&1.0 \\
\hline
\end{array}$$
} \caption{\label{VI}Performance of several families of valid
inequalities added to (F2)}
\end{table}

{
\def\arraystretch{0.5}
\def\tabcolsep{2pt}

\begin{table}
$$
{\small
\begin{array}{|c|rrr|rrr|rrr|}
\hline
 & \multicolumn{3}{|c|}{M=20} & \multicolumn{3}{|c|}{M=30} & \multicolumn{3}{|c|}{M=40} \\
\hline
 & p & \multicolumn{1}{c}{\bar{t}} & \multicolumn{1}{c|}{n} & \multicolumn{1}{c}{p} & \multicolumn{1}{c}{\bar{t}} &
\multicolumn{1}{c|}{n} & \multicolumn{1}{c}{p} & \multicolumn{1}{c}{\bar{t}} & \multicolumn{1}{c|}{n} \\
\hline
 \multicolumn{10}{|c|}{\hbox{Customers prefer closer sites and
self-service is allowed}} \\
\hline
\hbox{(F2)}                          & 2 & 10.7 & 420.2 & 3    & 63.4 & 948.2 &  2  & 379.9 & 14480.0 \\
\hbox{(F2)}(\ref{p101})              &   & 4.6  & 96.2  &      & 10.0 & 159.0 &     & 176.9 & 838.2  \\
\hbox{(F2)}(\ref{p101})(\ref{p104})  &   &  9.0 & 67.8  &      & 16.7 & 110.6 &     & 101.6  & 602.2 \\
\hline
\hbox{(F2)}                          & 3 & 7.3  &  317.8  & 6 & 4.5  & 113.4  &  4  & 265.8  & 2142.6 \\
\hbox{(F2)}(\ref{p101})              &   & 2.9  &  71.0        &  &3.0   &8.7          &  &51.6  &413.8  \\
\hbox{(F2)}(\ref{p101})(\ref{p104})  &   &  4.5 &  55.8   &   & 5.6  & 53.0   &     & 35.8   & 237.0 \\
\hline
\hbox{(F2)}   & 5 &  2.6 & 63.0   &10 & 1.3  & 35.4   &  8  & 16.6   & 207.8  \\
\hbox{(F2)}(\ref{p101})              &  &1.2&42.2         &  &1.1   &15.0          &  &13.1&143.8   \\
\hbox{(F2)}(\ref{p101})(\ref{p104})  &   &  1.5 & 17.4   &   & 2.7  & 13.0   &     & 19.2   & 119.8  \\
\hline
\hbox{(F2)}   & 7 &  1.4 & 27.4   &12 & 0.7  & 7.8    &  10 & 13.5   & 80.6   \\
\hbox{(F2)}(\ref{p101})              &  &0.4&10.0         &  &0.7   &2.3          &  &7.2   &86.5   \\
\hbox{(F2)}(\ref{p101})(\ref{p104})  &   &  0.8 & 8.6    &   & 1.4  & 4.2    &     & 10.5   & 56.2   \\
\hline
\hbox{(F2)}   & 10&  0.2 & 3.0    &15 & 1.9  & 13.0   & 16  & 2.0    & 12.2  \\
\hbox{(F2)}(\ref{p101})              &  &0.2 &5.0         &  &1.4   &11.4          &  &2.3   &18.5   \\
\hbox{(F2)}(\ref{p101})(\ref{p104})  &   &  0.6 & 3.8    &   & 1.2  & 6.2    &     & 5.0    & 15.8  \\
\hline
\hbox{(F2)}   & 12&  0.1 & 1.0    &22 & 0.2  & 1.0    & 20  & 3.8    & 1.8   \\
\hbox{(F2)}(\ref{p101})              &  &0.1 &1.0         &  &0.2   &1.0           &  &1.0   &1.0   \\
\hbox{(F2)}(\ref{p101})(\ref{p104})  &   &  0.3 & 1.0    &   & 0.5  & 1.0    &     & 3.1    & 1.4   \\
\hline
  \multicolumn{10}{|c|}{\hbox{Customers prefer closer sites but
self-service is forbidden}} \\
\hline
\hbox{(F2)}   & 2  &  28.2 &  2117.8 & 3  & 814.2 &  21899.4 & 2  & 2717.5 & 17868.5 \\
\hbox{(F2)}(\ref{p101})              &  &4.4 &193.4      &  &28.2 &825.4     &  &125.2 &876.6    \\
\hbox{(F2)}(\ref{p101})(\ref{p104})  &    &  4.0  &  140.6  &    & 29.9   & 665.4   &    & 108.0  & 719.0 \\
\hline
\hbox{(F2)}   & 3  &  16.6 &  1791.0 & 6  & 52.0  &  2749.8  & 4  & >1H    &  \multicolumn{1}{c|}{-}   \\
\hbox{(F2)}(\ref{p101})              &  &3.8 &275.4      &  &6.7  &304.6     &  &114.6 &1815.0   \\
\hbox{(F2)}(\ref{p101})(\ref{p104})  &    &  3.6  &  154.2  &    & 9.3   & 187.4   &    & 137.8 & 1382.7 \\
\hline
\hbox{(F2)}   & 5  &  2.5  &  364.6  & 10 & 4.5   &  420.4   & 8  & 846.7  & 28928.2 \\
\hbox{(F2)}(\ref{p101})              &  &3.8 &275.4      &  &3.1  &171.6     &  &23.1  &632.2    \\
\hbox{(F2)}(\ref{p101})(\ref{p104})  &    &  2.9  &  103.8  &    & 6.6   &  173.0   &    & 88.6  & 479.2 \\
\hline
\hbox{(F2)}   & 7  &  1.0  &  101.4  & 12 & 1.7   &  97.4    & 10 & 139.6  & 4933.4 \\
\hbox{(F2)}(\ref{p101})              &  &0.7 &67.4       &  &1.7  &72.2      &  &21.7  &740.6   \\
\hbox{(F2)}(\ref{p101})(\ref{p104})  &    &  1.5  &  45.4   &    & 3.7   &  65.0    &    & 15.7  & 231.7 \\
\hline
\hbox{(F2)}   & 10 &  0.4  &  28.6   & 15 & 2.2   &  21.8    & 16 & 4.2    & 136.0  \\
\hbox{(F2)}(\ref{p101})              &  &1.3 &17.4       &  &1.5  &15.0      &  &4.6   &75.0     \\
\hbox{(F2)}(\ref{p101})(\ref{p104})  &    &  0.8  &  20.2   &    & 1.4   &  13.0    &    & 7.4   & 70.6 \\
\hline
\hbox{(F2)}   & 12 &  0.4  &  10.6   & 22 & 0.2   &  1.8     & 20 & 2.0    & 53.8  \\
\hbox{(F2)}(\ref{p101})              &  &0.6 &6.6        &  &0.1  &1.0       &  &2.6   &19.0     \\
\hbox{(F2)}(\ref{p101})(\ref{p104})  &    &  0.9  &  4.6    &    & 1.3   &  1.0     &    & 4.1   & 12.0 \\
\hline
  \multicolumn{10}{|c|}{\hbox{Random preferences}} \\
\hline
\hbox{(F2)}   & 2  &  62.7 &  4650.6 & 3  & >1H   &    \multicolumn{1}{c|}{-}    & 2  & >1H    &  \multicolumn{1}{c|}{-}     \\
\hbox{(F2)}(\ref{p101})              &  &6.2 &372.6      &  &93.5 &2489.8      &  &>1H  &    \multicolumn{1}{c|}{-}     \\
\hbox{(F2)}(\ref{p101})(\ref{p104})  &    &  6.7  &  245.5  &    & 127.3 & 3204.6  &    & >1H   &   \multicolumn{1}{c|}{-}     \\
\hline
\hbox{(F2)}   & 3  &  15.5 &  1262.6 & 6  & 113.5 & 5463.8  & 4  & >1H    &   \multicolumn{1}{c|}{-}    \\
\hbox{(F2)}(\ref{p101})              &  &3.1 &287.4      &  &6.6  &245.4       &  &1342.9&16320.0  \\
\hbox{(F2)}(\ref{p101})(\ref{p104})  &    &  5.8  &  179.4  &    & 7.8   & 180.2   &    & 1253.0& 16697.7 \\
\hline
\hbox{(F2)}   & 5  &  2.4  &  329.4  & 10 & 2.7   & 137.8   & 8  & 212.8  & 4813.8  \\
\hbox{(F2)}(\ref{p101})              &  &1.2 &83.0       &  &2.2  &64.6        &  &32.3  &182.2    \\
\hbox{(F2)}(\ref{p101})(\ref{p104})  &    &  2.0  &  75.4   &    & 2.9   & 28.2    &    & 21.5  & 269.0   \\
\hline
\hbox{(F2)}   & 7  &  1.1  &  44.6   & 12 & 1.1   & 39.4    & 10 & 28.4   & 645.4 \\
\hbox{(F2)}(\ref{p101})              &  &2.1 &36.2       &  &1.4  &23.4        &  &10.3  &166.0    \\
\hbox{(F2)}(\ref{p101})(\ref{p104})  &    &  1.3  &  31.0   &    & 2.6   & 20.2    &    & 11.2  & 119.0 \\
\hline
\hbox{(F2)}   & 10 &  0.5  &  16.2   & 15 & 0.8   & 19.0    & 16 & 2.7    & 31.8  \\
\hbox{(F2)}(\ref{p101})              &  &1.0 &11.0       &  &1.1  &12.0        &  &32.0  &4.7     \\
\hbox{(F2)}(\ref{p101})(\ref{p104})  &    &  0.7  &  9.0    &    & 1.9   & 18.0    &    & 5.9   & 15.0  \\
\hline
\hbox{(F2)}   & 12 &  0.3  &  6.2    & 22 & 0.1   & 1.0     & 20 & 2.1    & 39.4   \\
\hbox{(F2)}(\ref{p101})              &  &0.3 &3.4        &  &0.1  &1.0         &  &2.8   &16.0     \\
\hbox{(F2)}(\ref{p101})(\ref{p104})  &    &  0.8  &  2.6    &    & 0.3   & 1.5     &    & 3.6   & 21.0   \\
\hline
\end{array}
}
$$
\caption{\label{KK} Comparison between (F2) with and without
separation of constraints (\ref{p104})}
\end{table}
}

\begin{table}
$$
{\small
\begin{array}{|c|rrrrr|rrrrr|rrrrr|}
\hline
 & \multicolumn{5}{|c|}{M=20} & \multicolumn{5}{|c|}{M=30} & \multicolumn{5}{|c|}{M=40} \\
\hline
 & p & \multicolumn{1}{c}{\bar{t}_P} & \multicolumn{1}{c}{\%{v}_0}&\multicolumn{1}{c}{\%{v}_1}&\multicolumn{1}{c|}{\bar{t}}
 & \multicolumn{1}{c}{p} &\multicolumn{1}{c}{\bar{t}_P} & \multicolumn{1}{c}{\%{v}_0}&\multicolumn{1}{c}{\%{v}_1}&\multicolumn{1}{c|}{\bar{t}}
 &\multicolumn{1}{c}{p}  &\multicolumn{1}{c}{\bar{t}_P} & \multicolumn{1}{c}{\%{v}_0}&\multicolumn{1}{c}{\%{v}_1}&\multicolumn{1}{c|}{\bar{t}}\\
\hline
  \multicolumn{16}{|c|}{\hbox{Customers prefer closer sites and
self-service is allowed}} \\
\hline
\hbox{(F2)}               &2 &\multicolumn{1}{c}{-}&\multicolumn{1}{c}{-}&\multicolumn{1}{c}{-}&10.8     &3 &\multicolumn{1}{c}{-}&\multicolumn{1}{c}{-}&\multicolumn{1}{c}{-}&63.4      &2 &\multicolumn{1}{c}{-}&\multicolumn{1}{c}{-}&\multicolumn{1}{c}{-}&2379.9\\
\hbox{(F2)(\ref{p101})}   &  &\multicolumn{1}{c}{-}&\multicolumn{1}{c}{-}&\multicolumn{1}{c}{-}&4.6      &  &\multicolumn{1}{c}{-}&\multicolumn{1}{c}{-}&\multicolumn{1}{c}{-}&10.0      &  &\multicolumn{1}{c}{-}&\multicolumn{1}{c}{-}&\multicolumn{1}{c}{-}&176.9 \\
\hbox{(F2)(\ref{p101})Pre}&  &0.4                  &58.6                 &15.9                 &1.5      &  &0.7                  &71.5                 &10.5                 &5.2       &  &0.9                  &78.1                 &14.3                 &59.0  \\
\hline
\hbox{(F2)}               &3 &\multicolumn{1}{c}{-}&\multicolumn{1}{c}{-}&\multicolumn{1}{c}{-}&7.3      &6 &\multicolumn{1}{c}{-}&\multicolumn{1}{c}{-}&\multicolumn{1}{c}{-}&4.5       &4 &\multicolumn{1}{c}{-}&\multicolumn{1}{c}{-}&\multicolumn{1}{c}{-}&265.8 \\
\hbox{(F2)(\ref{p101})}   &  &\multicolumn{1}{c}{-}&\multicolumn{1}{c}{-}&\multicolumn{1}{c}{-}&3.0      &  &\multicolumn{1}{c}{-}&\multicolumn{1}{c}{-}&\multicolumn{1}{c}{-}&3.0       &  &\multicolumn{1}{c}{-}&\multicolumn{1}{c}{-}&\multicolumn{1}{c}{-}&51.6  \\
\hbox{(F2)(\ref{p101})Pre}&  &0.5                  &59.1                 &11.2                 &1.2      &  &0.7                  &72.5                 &5.5                  &1.7       &  &1.0                  &78.4                 &8.2                  &16.6  \\
\hline
\hbox{(F2)}               &5 &\multicolumn{1}{c}{-}&\multicolumn{1}{c}{-}&\multicolumn{1}{c}{-}&2.6      &10&\multicolumn{1}{c}{-}&\multicolumn{1}{c}{-}&\multicolumn{1}{c}{-}&1.3       &8 &\multicolumn{1}{c}{-}&\multicolumn{1}{c}{-}&\multicolumn{1}{c}{-}&16.6  \\
\hbox{(F2)(\ref{p101})}   &  &\multicolumn{1}{c}{-}&\multicolumn{1}{c}{-}&\multicolumn{1}{c}{-}&1.2      &  &\multicolumn{1}{c}{-}&\multicolumn{1}{c}{-}&\multicolumn{1}{c}{-}&1.1       &  &\multicolumn{1}{c}{-}&\multicolumn{1}{c}{-}&\multicolumn{1}{c}{-}&13.1  \\
\hbox{(F2)(\ref{p101})Pre}&  &0.4                  &61.4                 &6.9                  &0.6      &  &0.8                  &72.5                 &3.9                  &0.9       &  &1.2                  &78.9                 &4.3                  &4.6   \\
\hline
\hbox{(F2)}               &7 &\multicolumn{1}{c}{-}&\multicolumn{1}{c}{-}&\multicolumn{1}{c}{-}&1.4      &12&\multicolumn{1}{c}{-}&\multicolumn{1}{c}{-}&\multicolumn{1}{c}{-}&0.7       &10&\multicolumn{1}{c}{-}&\multicolumn{1}{c}{-}&\multicolumn{1}{c}{-}&13.5  \\
\hbox{(F2)(\ref{p101})}   &  &\multicolumn{1}{c}{-}&\multicolumn{1}{c}{-}&\multicolumn{1}{c}{-}&0.4      &  &\multicolumn{1}{c}{-}&\multicolumn{1}{c}{-}&\multicolumn{1}{c}{-}&0.7       &  &\multicolumn{1}{c}{-}&\multicolumn{1}{c}{-}&\multicolumn{1}{c}{-}&7.2   \\
\hbox{(F2)(\ref{p101})Pre}&  &0.4                  &61.7                 &5.5                 &0.5       &  &0.8                  &73.8                 &3.5                  &0.8       &  &1.3                  &79.1                 &3.6                  &2.9   \\
\hline
\hbox{(F2)}               &10&\multicolumn{1}{c}{-}&\multicolumn{1}{c}{-}&\multicolumn{1}{c}{-}&0.2      &15&\multicolumn{1}{c}{-}&\multicolumn{1}{c}{-}&\multicolumn{1}{c}{-}&1.9       &16&\multicolumn{1}{c}{-}&\multicolumn{1}{c}{-}&\multicolumn{1}{c}{-}&2.0   \\
\hbox{(F2)(\ref{p101})}   &  &\multicolumn{1}{c}{-}&\multicolumn{1}{c}{-}&\multicolumn{1}{c}{-}&0.2      &  &\multicolumn{1}{c}{-}&\multicolumn{1}{c}{-}&\multicolumn{1}{c}{-}&1.4       &  &\multicolumn{1}{c}{-}&\multicolumn{1}{c}{-}&\multicolumn{1}{c}{-}&2.3   \\
\hbox{(F2)(\ref{p101})Pre}&  &0.4                  &62.8                 &5.0                  &0.4      &  &0.7                  &73.7                 &3.3                  &0.8       &  &1.1                  &79.4                 &2.6                  &1.2   \\
\hline
\hbox{(F2)}               &12&\multicolumn{1}{c}{-}&\multicolumn{1}{c}{-}&\multicolumn{1}{c}{-}&0.1      &22&\multicolumn{1}{c}{-}&\multicolumn{1}{c}{-}&\multicolumn{1}{c}{-}&0.2       &20&\multicolumn{1}{c}{-}&\multicolumn{1}{c}{-}&\multicolumn{1}{c}{-}&3.8   \\
\hbox{(F2)(\ref{p101})}   &  &\multicolumn{1}{c}{-}&\multicolumn{1}{c}{-}&\multicolumn{1}{c}{-}&0.1      &  &\multicolumn{1}{c}{-}&\multicolumn{1}{c}{-}&\multicolumn{1}{c}{-}&0.2       &  &\multicolumn{1}{c}{-}&\multicolumn{1}{c}{-}&\multicolumn{1}{c}{-}&1.0   \\
\hbox{(F2)(\ref{p101})Pre}&  & 2.2                 &70.5                 &5.0                  &2.5      &  & 1.6                 &79.6                 &3.3                  &1.7       &  &1.7                  &79.9                 &3.3                  &2.0   \\
\hline
  \multicolumn{16}{|c|}{\hbox{Customers prefer closer sites but
self-service is forbidden}} \\
\hline
\hbox{(F2)}               &2 &\multicolumn{1}{c}{-}&\multicolumn{1}{c}{-}&\multicolumn{1}{c}{-}&28.1     &3 &\multicolumn{1}{c}{-}&\multicolumn{1}{c}{-}&\multicolumn{1}{c}{-}&814.1     &2 &\multicolumn{1}{c}{-}&\multicolumn{1}{c}{-}&\multicolumn{1}{c}{-}&2717.5\\
\hbox{(F2)(\ref{p101})}   &  &\multicolumn{1}{c}{-}&\multicolumn{1}{c}{-}&\multicolumn{1}{c}{-}&4.4      &  &\multicolumn{1}{c}{-}&\multicolumn{1}{c}{-}&\multicolumn{1}{c}{-}&28.1      &  &\multicolumn{1}{c}{-}&\multicolumn{1}{c}{-}&\multicolumn{1}{c}{-}&125.2  \\
\hbox{(F2)(\ref{p101})Pre}&  &0.5                  &57.5                 &14.4                 &2.7      &  &0.8                  &70.9                 &9.8                  &20.5      &  &1.0                  &77.8                 &13.6                 &74.6  \\
\hline
\hbox{(F2)}               &3 &\multicolumn{1}{c}{-}&\multicolumn{1}{c}{-}&\multicolumn{1}{c}{-}&16.6     &6 &\multicolumn{1}{c}{-}&\multicolumn{1}{c}{-}&\multicolumn{1}{c}{-}&52.0      &4 &\multicolumn{1}{c}{-}&\multicolumn{1}{c}{-}&\multicolumn{1}{c}{-}&>1H  \\
\hbox{(F2)(\ref{p101})}   &  &\multicolumn{1}{c}{-}&\multicolumn{1}{c}{-}&\multicolumn{1}{c}{-}&3.8      &  &\multicolumn{1}{c}{-}&\multicolumn{1}{c}{-}&\multicolumn{1}{c}{-}&6.7       &  &\multicolumn{1}{c}{-}&\multicolumn{1}{c}{-}&\multicolumn{1}{c}{-}&114.6  \\
\hbox{(F2)(\ref{p101})Pre}&  &0.5                  &58.2                 &9.8                  &2.4      &  &1.0                  &71.3                 &4.4                  &4.2       &  &1.3                  &77.8                 &7.7                  &90.8  \\
\hline
\hbox{(F2)}               &5 &\multicolumn{1}{c}{-}&\multicolumn{1}{c}{-}&\multicolumn{1}{c}{-}&2.5      &10&\multicolumn{1}{c}{-}&\multicolumn{1}{c}{-}&\multicolumn{1}{c}{-}&4.5       &8 &\multicolumn{1}{c}{-}&\multicolumn{1}{c}{-}&\multicolumn{1}{c}{-}&846.7 \\
\hbox{(F2)(\ref{p101})}   &  &\multicolumn{1}{c}{-}&\multicolumn{1}{c}{-}&\multicolumn{1}{c}{-}&3.8      &  &\multicolumn{1}{c}{-}&\multicolumn{1}{c}{-}&\multicolumn{1}{c}{-}&3.1       &  &\multicolumn{1}{c}{-}&\multicolumn{1}{c}{-}&\multicolumn{1}{c}{-}&23.1  \\
\hbox{(F2)(\ref{p101})Pre}&  &0.5                  &59.5                 &5.4                  &1.2      &  &0.9                  &71.3                 &2.6                  &1.8       &  &1.7                  &78.0                 &3.5                  &22.1  \\
\hline
\hbox{(F2)}               &7 &\multicolumn{1}{c}{-}&\multicolumn{1}{c}{-}&\multicolumn{1}{c}{-}&1.0      &12&\multicolumn{1}{c}{-}&\multicolumn{1}{c}{-}&\multicolumn{1}{c}{-}&1.7       &10&\multicolumn{1}{c}{-}&\multicolumn{1}{c}{-}&\multicolumn{1}{c}{-}&139.6 \\
\hbox{(F2)(\ref{p101})}   &  &\multicolumn{1}{c}{-}&\multicolumn{1}{c}{-}&\multicolumn{1}{c}{-}&0.7      &  &\multicolumn{1}{c}{-}&\multicolumn{1}{c}{-}&\multicolumn{1}{c}{-}&1.7       &  &\multicolumn{1}{c}{-}&\multicolumn{1}{c}{-}&\multicolumn{1}{c}{-}&21.7  \\
\hbox{(F2)(\ref{p101})Pre}&  &0.5                  &59.3                 &3.4                  &0.7      &  &0.9                  &71.9                 &1.8                  &1.2       &  &1.9                  &78.1                 &2.7                  &14.3  \\
\hline
\hbox{(F2)}               &10&\multicolumn{1}{c}{-}&\multicolumn{1}{c}{-}&\multicolumn{1}{c}{-}&0.4      &15&\multicolumn{1}{c}{-}&\multicolumn{1}{c}{-}&\multicolumn{1}{c}{-}&2.2       &16&\multicolumn{1}{c}{-}&\multicolumn{1}{c}{-}&\multicolumn{1}{c}{-}&4.2   \\
\hbox{(F2)(\ref{p101})}   &  &\multicolumn{1}{c}{-}&\multicolumn{1}{c}{-}&\multicolumn{1}{c}{-}&1.3      &  &\multicolumn{1}{c}{-}&\multicolumn{1}{c}{-}&\multicolumn{1}{c}{-}&1.5       &  &\multicolumn{1}{c}{-}&\multicolumn{1}{c}{-}&\multicolumn{1}{c}{-}&4.6   \\
\hbox{(F2)(\ref{p101})Pre}&  &0.5                  &60.9                 &2.0                  &0.5      &  &0.8                  &71.6                 &1.4                  &1.0       &  &1.5                  &78.0                 &1.5                  &2.1   \\
\hline
\hbox{(F2)}               &12&\multicolumn{1}{c}{-}&\multicolumn{1}{c}{-}&\multicolumn{1}{c}{-}&0.4      &22&\multicolumn{1}{c}{-}&\multicolumn{1}{c}{-}&\multicolumn{1}{c}{-}&0.2       &20&\multicolumn{1}{c}{-}&\multicolumn{1}{c}{-}&\multicolumn{1}{c}{-}&2.0   \\
\hbox{(F2)(\ref{p101})}   &  &\multicolumn{1}{c}{-}&\multicolumn{1}{c}{-}&\multicolumn{1}{c}{-}&0.6      &  &\multicolumn{1}{c}{-}&\multicolumn{1}{c}{-}&\multicolumn{1}{c}{-}&0.1       &  &\multicolumn{1}{c}{-}&\multicolumn{1}{c}{-}&\multicolumn{1}{c}{-}&2.6   \\
\hbox{(F2)(\ref{p101})Pre}&  &2.0                  &67.4                 &1.1                  &2.1      &  &1.5                  &76.0                 &0.0                  &1.7       &  &1.7                  &78.1                 &1.1                 &2.0   \\
\hline
 \multicolumn{16}{|c|}{\hbox{Random preferences}} \\
\hline
\hbox{(F2)}               &2 &\multicolumn{1}{c}{-}&\multicolumn{1}{c}{-}&\multicolumn{1}{c}{-}&62.7     &3 &\multicolumn{1}{c}{-}&\multicolumn{1}{c}{-}&\multicolumn{1}{c}{-}&>1H       &2 &\multicolumn{1}{c}{-}&\multicolumn{1}{c}{-}&\multicolumn{1}{c}{-}&>1H    \\
\hbox{(F2)(\ref{p101})}   &  &\multicolumn{1}{c}{-}&\multicolumn{1}{c}{-}&\multicolumn{1}{c}{-}&6.2      &  &\multicolumn{1}{c}{-}&\multicolumn{1}{c}{-}&\multicolumn{1}{c}{-}&93.5      &  &\multicolumn{1}{c}{-}&\multicolumn{1}{c}{-}&\multicolumn{1}{c}{-}&>1H    \\
\hbox{(F2)(\ref{p101})Pre}&  &1.0                  &65.8                 &14.6                 &6.1      &  &1.8                  &72.0                 &10.6                 &65.2      &  &2.0                  &78.2                 &13.7                 &2005.5 \\
\hline
\hbox{(F2)}               &3 &\multicolumn{1}{c}{-}&\multicolumn{1}{c}{-}&\multicolumn{1}{c}{-}&15.5     &6 &\multicolumn{1}{c}{-}&\multicolumn{1}{c}{-}&\multicolumn{1}{c}{-}&113.5     &4 &\multicolumn{1}{c}{-}&\multicolumn{1}{c}{-}&\multicolumn{1}{c}{-}&>1H     \\
\hbox{(F2)(\ref{p101})}   &  &\multicolumn{1}{c}{-}&\multicolumn{1}{c}{-}&\multicolumn{1}{c}{-}&3.1      &  &\multicolumn{1}{c}{-}&\multicolumn{1}{c}{-}&\multicolumn{1}{c}{-}&6.6       &  &\multicolumn{1}{c}{-}&\multicolumn{1}{c}{-}&\multicolumn{1}{c}{-}&1342.9   \\
\hbox{(F2)(\ref{p101})Pre}&  &1.1                  &66.6                 &9.0                  &4.3      &  &3.1                  &74.5                 &4.1                  &6.9       &  &4.3                  &79.6                 &7.8                  &397.5 \\
\hline
\hbox{(F2)}               &5 &\multicolumn{1}{c}{-}&\multicolumn{1}{c}{-}&\multicolumn{1}{c}{-}&2.4      &10&\multicolumn{1}{c}{-}&\multicolumn{1}{c}{-}&\multicolumn{1}{c}{-}&2.7       &8 &\multicolumn{1}{c}{-}&\multicolumn{1}{c}{-}&\multicolumn{1}{c}{-}&212.8    \\
\hbox{(F2)(\ref{p101})}   &  &\multicolumn{1}{c}{-}&\multicolumn{1}{c}{-}&\multicolumn{1}{c}{-}&1.2      &  &\multicolumn{1}{c}{-}&\multicolumn{1}{c}{-}&\multicolumn{1}{c}{-}&2.2       &  &\multicolumn{1}{c}{-}&\multicolumn{1}{c}{-}&\multicolumn{1}{c}{-}&32.3\\
\hbox{(F2)(\ref{p101})Pre}&  &1.0                  &68.5                 &4.3                  &1.7      &  &3.4                  &77.2                 &1.5                  &4.0       &  &11.3                 &80.8                 &3.1                  &20.7\\
\hline
\hbox{(F2)}               &7 &\multicolumn{1}{c}{-}&\multicolumn{1}{c}{-}&\multicolumn{1}{c}{-}&1.1      &12&\multicolumn{1}{c}{-}&\multicolumn{1}{c}{-}&\multicolumn{1}{c}{-}&1.1       &10&\multicolumn{1}{c}{-}&\multicolumn{1}{c}{-}&\multicolumn{1}{c}{-}&28.4     \\
\hbox{(F2)(\ref{p101})}   &  &\multicolumn{1}{c}{-}&\multicolumn{1}{c}{-}&\multicolumn{1}{c}{-}&2.1      &  &\multicolumn{1}{c}{-}&\multicolumn{1}{c}{-}&\multicolumn{1}{c}{-}&1.4       &  &\multicolumn{1}{c}{-}&\multicolumn{1}{c}{-}&\multicolumn{1}{c}{-}&10.3  \\
\hbox{(F2)(\ref{p101})Pre}&  &1.0                  &69.8                 &2.4                  &1.3      &  &3.2                  &77.8                 &1.2                  &3.6       &  &16.8                 &82.1                 &2.1                  &19.3  \\
\hline
\hbox{(F2)}               &10&\multicolumn{1}{c}{-}&\multicolumn{1}{c}{-}&\multicolumn{1}{c}{-}&0.5      &15&\multicolumn{1}{c}{-}&\multicolumn{1}{c}{-}&\multicolumn{1}{c}{-}&0.8       &16&\multicolumn{1}{c}{-}&\multicolumn{1}{c}{-}&\multicolumn{1}{c}{-}&2.7      \\
\hbox{(F2)(\ref{p101})}   &  &\multicolumn{1}{c}{-}&\multicolumn{1}{c}{-}&\multicolumn{1}{c}{-}&1.0      &  &\multicolumn{1}{c}{-}&\multicolumn{1}{c}{-}&\multicolumn{1}{c}{-}&1.1       &  &\multicolumn{1}{c}{-}&\multicolumn{1}{c}{-}&\multicolumn{1}{c}{-}&32.0  \\
\hbox{(F2)(\ref{p101})Pre}&  &0.9                  &71.8                 &1.2                  &0.9      &  &2.5                  &78.7                 &0.6                  &2.6       &  &17.2                 &83.4                 &1.1                  &17.7 \\
\hline
\hbox{(F2)}               &12&\multicolumn{1}{c}{-}&\multicolumn{1}{c}{-}&\multicolumn{1}{c}{-}&0.3      &22&\multicolumn{1}{c}{-}&\multicolumn{1}{c}{-}&\multicolumn{1}{c}{-}&0.1       &20&\multicolumn{1}{c}{-}&\multicolumn{1}{c}{-}&\multicolumn{1}{c}{-}&2.1      \\
\hbox{(F2)(\ref{p101})}   &  &\multicolumn{1}{c}{-}&\multicolumn{1}{c}{-}&\multicolumn{1}{c}{-}&0.3      &  &\multicolumn{1}{c}{-}&\multicolumn{1}{c}{-}&\multicolumn{1}{c}{-}&0.1       &  &\multicolumn{1}{c}{-}&\multicolumn{1}{c}{-}&\multicolumn{1}{c}{-}&2.8   \\
\hbox{(F2)(\ref{p101})Pre}&  &1.6                  &73.2                &0.5                  &1.8      &  &1.5                  &81.7                 &0.0                  &1.8       &  &11.5                 &84.4                 &0.9                  &11.8  \\
\hline
\end{array}
}
$$
\caption{\label{f2Pre} Comparison between (F2) with and without
preprocessing}
\end{table}

In Table \ref{f3mejor}, the results obtained for formulation (F3)
with and without the modified objective function have been compared.

{\small
\begin{table}
$$
\hspace{-2cm}
\begin{array}{|c|rrrrr|rrrrr|rrrrr|}
\hline
\multicolumn{1}{|c|}{} & \multicolumn{5}{|c|}{M=20} & \multicolumn{5}{|c|}{M=30} & \multicolumn{5}{|c|}{M=40} \\
\hline
 & \multicolumn{1}{c}{p} & \multicolumn{1}{c}{LP} & \multicolumn{1}{c}{\bar{t}} & \multicolumn{1}{c}{\sigma_t} &
\multicolumn{1}{c|}{n} & \multicolumn{1}{c}{p} &
\multicolumn{1}{c}{LP} & \multicolumn{1}{c}{\bar{t}} &
\multicolumn{1}{c}{\sigma_t} & \multicolumn{1}{c|}{n} &
\multicolumn{1}{c}{p} & \multicolumn{1}{c}{LP} &
\multicolumn{1}{c}{\bar{t}} & \multicolumn{1}{c}{\sigma_t} & \multicolumn{1}{c|}{n} \\
\hline
  \multicolumn{16}{|c|}{\hbox{Customers prefer closer sites and
self-service is allowed}} \\
\hline
\hbox{(F3)}   &2 &0.0 &1.5 &0.1 &177.0  &       3&0.1 &22.7  &3.6  &718.2    &       2 &0.5 &50.7  &4.5   &1077.0 \\
\hbox{(F3R)}  &  &0.0 &1.4 &0.2 &153.4  &        &0.1 &19.2  &3.4  &617.0    &         &0.2 &42.3  &4.3   &577.8  \\
\hline
\hbox{(F3)}   &3 &0.0 &2.5 &0.3 &327.8  &       6&0.1 &81.9  &10.8 &4357.8   &       4 &0.2 &263.6 &29.7  &3373.8  \\
\hbox{(F3R)}  &  &0.0 &2.2 &0.4 &291.4  &        &0.1 &65.8  &25.5 &3955.0   &         &0.2 &218.1 &32.5  &3568.6\\
\hline
\hbox{(F3)}   &5 &0.0 &4.4&1.3  &917.0  &      10&0.1 &174.5 &84.1 &17987.0  &       8 &0.2 &2443.9&357.9 &57277.8\\
\hbox{(F3R)}  &  &0.0 &3.6&1.1  &718.6  &        &0.1 &141.5 &13.3 &10941.0  &         &0.1 &1568.4&165.8 &37586.2\\
\hline
\hbox{(F3)}   &7 &0.0 &4.4&0.8  &1255.8 &      12&0.1 &425.6 &260.4&51375.0  &       10&0.2 &>1H  &   -    &  - \\
\hbox{(F3R)}  &  &0.0 &4.0&0.7  &874.6  &        &0.0 &270.5 &161.4&23519.8  &         &0.1 &3200.6&626.6 &81377.5\\
\hline
\hbox{(F3)}   &10&0.0 &9.5 &2.0 &5579.8 &      15&0.1 &996.8 &25.0 &249648.0 &       16&0.2 &>1H  &   -    &   - \\
\hbox{(F3R)}  &  &0.0 &7.6 &0.9 &2937.0 &        &0.0 &686.0 &123.2&107674.0 &          &0.2 &>1H  &  -     &   -   \\
\hline
\hbox{(F3)}   &12&0.0 &16.4&1.5 &14263.0&      22&0.1 &3298.1&150.8&1999284.0&       20&0.2 &>1H  &    -   &    - \\
\hbox{(F3R)}  &  &0.0 &17.8&2.4 &10633.4&        &\multicolumn{1}{c}{\hbox{OOM}} & -     & -    & -       &         &0.2 &>1H  &  -     &   -   \\
\hline
 \multicolumn{16}{|c|}{\hbox{Customers prefer closer sites but
self-service is forbidden}} \\
\hline
\hbox{(F3)}   &2 &0.0&1.4  &0.3  &167.8 &     3 &0.1 &26.5 &2.8  &861.9 &     2 &0.3 &53.1  &11.2  &482.2          \\
\hbox{(F3R)}  &  &0.0&1.4  &0.4  &137.4 &       &0.1 &22.2 &4.0  &786.2 &       &0.2 &39.3  &6.1   &409.4 \\
\hline
\hbox{(F3)}   &3 &0.0&2.5  &0.4  &315.4 &     6 &0.1 &66.9 &25.0 &3258.6&     4 &0.2 &326.3 &84.6  &4029.4         \\
\hbox{(F3R)}  &  &0.0&2.2  &0.5  &279.4 &       &0.1 &57.8 &29.9 &2917.8&       &0.2 &214.4 &67.7  &3193.8\\
\hline
\hbox{(F3)}   &5 &0.0&3.2  &0.5  &597.0 &     10&0.1 &82.3 &6.6  &6807.0&     8 &0.2 &2105.2&302.5 &38412.2       \\
\hbox{(F3R)}  &  &0.0&2.7  &0.1  &466.6  &      &0.1 &64.4 &7.1  &4516.0&       &0.1 &1281.3&360.6 &28625.8\\
\hline
\hbox{(F3)}   &7 &0.0&2.5  &1.1  &595.4 &     12&0.1 &48.8 &17.7 &4261.0&     10&0.2 &2763.5&449.3 &67043.0\\
\hbox{(F3R)}  &  &0.0&2.3  &0.6  &474.2 &       &0.0 &31.6 &8.5  &2208.2&       &0.1 &1829.3&571.6 &45384.0\\
\hline
\hbox{(F3)}   &10&0.0&1.3  &0.7  &411.0 &     15&0.1 &15.8 &0.6  &2054.0&     16&0.2 &1677.0&1147.9&43331.5       \\
\hbox{(F3R)}  &  &0.0&1.1  &0.7  &281.4 &       &0.0 &8.6  &0.6  &782.0 &       &0.1 &1346.5&948.2 &40464.6\\
\hline
\hbox{(F3)}   &12&0.0&0.4  &0.2  &107.4 &     22&0.0 &0.2  &0.0 &3.0    &     20&0.2 &127.5 &26.7  &6387.0           \\
\hbox{(F3R)}  &  &0.0&0.4  &0.2  &96.2  &       &0.0 &0.4  &0.1 &30.0     &       &0.1 &96.2  &32.0  &3906.0\\
\hline
  \multicolumn{16}{|c|}{\hbox{Random preferences}} \\
\hline
\hbox{(F3)}   &2 &0.0&2.0 &0.3  &242.6 &     3 &0.1 &34.0 &8.0  &1033.4&      2 &0.0&96.6  &13.0 &1150.2 \\
\hbox{(F3R)}  &  &0.0&1.7 &0.2  &175.4 &       &0.1 &23.9 &2.8  &699.4 &        &0.2&70.4  &4.5  &1028.2 \\
\hline
\hbox{(F3)}   &3 &0.0&2.5 &0.1  &293.4 &     6 &0.1 &63.2 &10.6 &2436.6&      4 &0.3&507.0 &57.1 &5660.2\\
\hbox{(F3R)}  &  &0.0&2.2 &0.3  &257.4 &       &0.1 &53.5 &15.7 &2183.8&        &0.2&319.4 &66.4 &4102.5\\
\hline
\hbox{(F3)}   &5 &0.0&2.3 &0.4  &346.2 &     10&0.1 &27.8&0.1   &1443.0&      8 &0.2&862.9 &402.6&12009.4\\
\hbox{(F3R)}  &  &0.0&2.1 &0.4  &317.8 &       &0.1 &29.8&8.4   &1504.0&        &0.2&447.3 &259.3&6487.8\\
\hline
\hbox{(F3)}   &7 &0.0&1.7 &0.6  &329.4 &     12&0.1 &20.5 &16.3 &1495.4&      10&0.2&651.4 &139.3&11759.0\\
\hbox{(F3R)}  &  &0.0&1.2 &0.5  &204.6 &       &0.0 &15.3 &10.0 &919.0 &        &0.2&423.3 &48.0 &7875.0 \\
\hline
\hbox{(F3)}   &10&0.0&0.6 &0.4  &166.2 &     15&0.1 &6.4  &7.1  &648.0 &      16&0.2&336.1 &202.0&12150.0\\
\hbox{(F3R)}  &  &0.0&0.7 &0.4  &160.6 &       &0.1 &3.8  &3.4  &326.0 &        &0.1&388.1 &244.2&10534.0\\
\hline
\hbox{(F3)}   &12&0.0&0.2 &0.2  &51.0  &     22&0.0 &0.2  &0.1  &3.0   &      20&0.2&71.1  &8.3  &3689.0\\
\hbox{(F3R)}  &  &0.0&0.2 &0.0  &41.4  &       &0.0 &0.3  &0.1  &30.0   &         &0.1&59.6  &41.5 &2340.0\\
\hline
\end{array}$$
\caption{\label{f3mejor} Comparison between (F3) with and without a
modification in the objective function}
\end{table}
}

A comparison between (F5.1) with and without preprocessing is shown
in Table \ref{f51Pre}.

\begin{table}
$$
{\small
\begin{array}{|c|rrr|rrr|rrr|}
\hline
 & \multicolumn{3}{|c|}{M=20} & \multicolumn{3}{|c|}{M=30} & \multicolumn{3}{|c|}{M=40} \\
\hline
 & p & \multicolumn{1}{c}{\bar{t}_P} &\multicolumn{1}{c|}{\bar{t}}
 & \multicolumn{1}{c}{p} &\multicolumn{1}{c}{\bar{t}_P} & \multicolumn{1}{c|}{\bar{t}}
 &\multicolumn{1}{c}{p}  &\multicolumn{1}{c}{\bar{t}_P} & \multicolumn{1}{c|}{\bar{t}}\\
\hline
  \multicolumn{10}{|c|}{\hbox{Customers prefer closer sites and
self-service is allowed}} \\
\hline
\hbox{(F5.1)}    &2 &\multicolumn{1}{c}{-}&1.1    &3 &\multicolumn{1}{c}{-}&25.8   &2 &\multicolumn{1}{c}{-}&51.2    \\
\hbox{(F5.1)Pre} &  &2.3                  &3.4    &  &2.4                  &38.2   &  &2.9                  &56.3    \\
\hline
\hbox{(F5.1)}    &3 &\multicolumn{1}{c}{-}&2.4    &6 &\multicolumn{1}{c}{-}&145.8  &4 &\multicolumn{1}{c}{-}&496.7   \\
\hbox{(F5.1)Pre} &  &2.0                  &4.0    &  &2.0                  &51.5   &  &2.8                  &881.6   \\
\hline
\hbox{(F5.1)}    &5 &\multicolumn{1}{c}{-}&4.9    &10&\multicolumn{1}{c}{-}&340.6  &8 &\multicolumn{1}{c}{-}&>1H     \\
\hbox{(F5.1)Pre} &  &1.9                  &3.2    &  &1.8                  &2.7    &  &2.3                  &>1H     \\
\hline
\hbox{(F5.1)}    &7 &\multicolumn{1}{c}{-}&6.1    &12&\multicolumn{1}{c}{-}&387.7  &10&\multicolumn{1}{c}{-}&>1H     \\
\hbox{(F5.1)Pre} &  &1.4                  &1.6    &  &1.8                  &52.2   &  &2.3                  &>1H     \\
\hline
\hbox{(F5.1)}    &10&\multicolumn{1}{c}{-}&18.0   &15&\multicolumn{1}{c}{-}&1701.9 &16&\multicolumn{1}{c}{-}&>1H     \\
\hbox{(F5.1)Pre} &  &1.5                  &1.6    &  &0.8                  &1058.9 &  &1.9                  &>1H     \\
\hline
\hbox{(F5.1)}    &12&\multicolumn{1}{c}{-}&33.2   &22&\multicolumn{1}{c}{-}&>1H    &20&\multicolumn{1}{c}{-}&>1H    \\
\hbox{(F5.1)Pre} &  &1.2                  &1.3    &  &1.3                  &1.6    &  &2.0                  &>1H     \\
\hline
 \multicolumn{10}{|c|}{\hbox{Customers prefer closer sites but
self-service is forbidden}} \\
\hline
\hbox{(F5.1)}    &2 &\multicolumn{1}{c}{-}&1.0    &3 &\multicolumn{1}{c}{-}&28.1   &2 &\multicolumn{1}{c}{-}&48.3  \\
\hbox{(F5.1)Pre} &  &2.1                  &4.1    &  &1.1                  &57.7   &  &1.6                  &81.4  \\
\hline
\hbox{(F5.1)}    &3 &\multicolumn{1}{c}{-}&2.4    &6 &\multicolumn{1}{c}{-}&92.0   &4 &\multicolumn{1}{c}{-}&521.6 \\
\hbox{(F5.1)Pre} &  &1.9                  &5.5    &  &1.0                  &176.0  &  &1.8                  &1779.1 \\
\hline
\hbox{(F5.1)}    &5 &\multicolumn{1}{c}{-}&4.3    &10&\multicolumn{1}{c}{-}&231.8  &8 &\multicolumn{1}{c}{-}&>1H   \\
\hbox{(F5.1)Pre} &  &1.7                  &6.0    &  &0.9                  &310.4  &  &2.4                  &>1H   \\
\hline
\hbox{(F5.1)}    &7 &\multicolumn{1}{c}{-}&4.7    &12&\multicolumn{1}{c}{-}&115.0  &10&\multicolumn{1}{c}{-}&>1H   \\
\hbox{(F5.1)Pre} &  &1.7                  &5.7    &  &1.0                  &81.2   &  &2.2                  &>1H   \\
\hline
\hbox{(F5.1)}    &10&\multicolumn{1}{c}{-}&1.5    &15&\multicolumn{1}{c}{-}&27.7   &16&\multicolumn{1}{c}{-}&>1H   \\
\hbox{(F5.1)Pre} &  &1.6                  &2.4    &  &0.9                  &17.5   &  &2.1                  &2549.5 \\
\hline
\hbox{(F5.1)}    &12&\multicolumn{1}{c}{-}&0.4    &22&\multicolumn{1}{c}{-}&0.1    &20&\multicolumn{1}{c}{-}&353.4 \\
\hbox{(F5.1)Pre} &  &1.2                  &1.3    &  &0.6                  &0.7    &  &1.5                  &9.9\\
\hline
 & \multicolumn{9}{c|}{\hbox{Random preferences}} \\
\hline
\hbox{(F5.1)}    &2 &\multicolumn{1}{c}{-}&1.5    &3 &\multicolumn{1}{c}{-}&42.7    &2 &\multicolumn{1}{c}{-}&100.3  \\
\hbox{(F5.1)Pre} &  &1.7                  &4.2    &  &1.3                  &106.3   &  &1.8                  &150.4 \\
\hline
\hbox{(F5.1)}    &3 &\multicolumn{1}{c}{-}&74.8   &6 &\multicolumn{1}{c}{-}&118.7   &4 &\multicolumn{1}{c}{-}&917.3  \\
\hbox{(F5.1)Pre} &  &1.8                  &5.9    &  &2.7                  &163.9   &  &3.1                  &>1H       \\
\hline
\hbox{(F5.1)}    &5 &\multicolumn{1}{c}{-}&2.8    &10&\multicolumn{1}{c}{-}&43.4    &8 &\multicolumn{1}{c}{-}&1969.2 \\
\hbox{(F5.1)Pre} &  &1.6                  &4.6    &  &3.6                   &56.1    &  &4.7                 &3066.1 \\
\hline
\hbox{(F5.1)}    &7 &\multicolumn{1}{c}{-}&1.8    &12&\multicolumn{1}{c}{-}&36.7    &10&\multicolumn{1}{c}{-}&2390.3 \\
\hbox{(F5.1)Pre} &  &1.5                  &2.9    &  &3.5                  &48.6    &  &8.8                  &1355.5 \\
\hline
\hbox{(F5.1)}    &10&\multicolumn{1}{c}{-}&0.6    &15&\multicolumn{1}{c}{-}&6.3     &16&\multicolumn{1}{c}{-}&798.8  \\
\hbox{(F5.1)Pre} &  &1.5                  &1.7    &  &4.2                  &4.8     &  &20.8                 &642.9   \\
\hline
\hbox{(F5.1)}    &12&\multicolumn{1}{c}{-}&0.2    &22&\multicolumn{1}{c}{-}&0.1     &20&\multicolumn{1}{c}{-}&108.8  \\
\hbox{(F5.1)Pre} &  &1.3                  &1.5    &  &2.0                  &2.1     &  &20.2                 &48.4   \\
\hline
\end{array}
}
$$
\caption{\label{f51Pre} Comparison between (F5.1) with and without
preprocessing}
\end{table}

\end{document}